\documentclass[12pt]{iopart}
\usepackage{amssymb}
\usepackage{diagram}
\usepackage[dvips]{rotating}
\usepackage{epsfig}

\newtheorem{defn}{Definition}
\newtheorem{thm}{Theorem}
\newtheorem{lem}{Lemma}

\newtheorem{prop}{Proposition}

\newcommand{\pf}{{\it Proof:}}
\newcommand{\ass}{\mathfrak{a}}
\newcommand{\iss}{\mathfrak{i}}
\newcommand{\qss}{\mathfrak{q}}
\newcommand{\ph}[1]{\phantom{#1}}
\newcommand{\com}{\mathfrak{c}}
\newcommand{\lid}{\mathfrak{l}}
\newcommand{\rid}{\mathfrak{r}}

\newcommand{\can}{\textrm{{\bf can}}}
\newcommand{\sym }{\textrm{{\bf sym}}}

\newcommand{\Endo}{\mbox{\textit{End}}}
\newcommand{\ubar}{\underline{\ph{S}}}
\newcommand{\maxm}{\textrm{max}}

\newcommand{\comb}[2]{{\tiny \left( \begin{array}{c} #1 \\ #2
\end{array}\right) }}
\newcommand{\Comb}[2]{\left( \begin{array}{c} #1 \\ #2
\end{array}\right) }
\newcommand{\hght}[1]{\underline{\overline{ #1 }}}
\newcommand{\nodes}[1]{|\hght{ #1 }|}
\newcommand{\BTree}{\textrm{\bf BTree}}
\newcommand{\IRBTree}{\textrm{\bf IRBTree}}
\newcommand{\IRNBTree}{\textrm{\bf IRNBTree}}
\newcommand{\RBTree}{\textrm{\bf RBTree}}
\newcommand{\RRBTree}{\textrm{\bf RRBTree}}
\newcommand{\RRNBTree}{\textrm{\bf RRNBTree}}

\newcommand{\NRBTree}{\textrm{\bf NRBTree}}
\newcommand{\NRRBTree}{\textrm{\bf NRBTree}}

\newcommand{\NRRNBTree}{\textrm{\bf NRRNBTree}}
\newcommand{\qB}{\textrm{\bf BRRBTree}}
\newcommand{\cB}{\mathcal{B}}
\newcommand{\cC}{\mathcal{C}}
\newcommand{\cD}{\mathcal{D}}
\newcommand{\cLO}{\mathcal{LO}}
\newcommand{\bN}{\mathbb{N}}
\newcommand{\COLbin}[2]{{\tiny \begin{array}{c} #1 \\ #2 \end{array}}}

\newcommand{\It}{\textrm{It}}

\begin{document}

\title{Natural Associativity without the Pentagon Condition}
\author{William P. Joyce}
\address{Department of Physics and Astronomy, University of Canterbury,
Private Bag 4800, Christchurch, New Zealand.\\
  email:w.joyce@phys.canterbury.ac.nz}

\begin{abstract}
A premonoidal category is equipped only with a bifunctor and a natural
isomorphism for associativity. We introduce a (deformation) natural
automorphism, $\qss $, representing the deviation from the Pentagon
condition. We uncover a binary tree representation for all diagrams
involving $\ass $ and $\qss $ and provide a link to permutations and
linear orderings. This leads to other notions of premonoidalness. We
define these notions and prove coherence results for each.
\end{abstract}

\pacno{02.10.Ws, 02.20.Qs, 02.20.Fh, 02.20.Df, 03.65.Fd, 31.15.Hz}

\maketitle

\section{Introduction}                                           %

An initial motivation for this paper is to address a simple
short--coming of monoidal categories. Namely the construction of a
purely fermionic statistics. Furthermore, to generalise statistical
structures in a physically meaningful way. In other words a
commutativity constraint given by $a\otimes b= -b\otimes a$. According
to the Hexagon diagram this requires associativity to be given by
$(a\otimes b)\otimes c=-a\otimes (b\otimes c)$. The Pentagon diagram
rules this possibility out. This paper describes the structure required
for such a choice and its natural extensions. We begin with a natural
isomorphism for associativity onto which we progressively add other
structures. This basic structure is called premonoidal since the
Pentagon diagram does not hold and the notion of identity is omitted.

The notion of coherence also changes. It is not possible to construct a
template premonoidal structure which can be shown to be related by a
premonoidal functor equivalence to every other premonoidal structure.
Instead we only ask for a groupoid whose diagrams somehow encode the
coherent diagrams of the premonoidal structure; and a functor for
interpreting these diagrams as diagrams in the category of interest.
Coherence asserts that all diagrams arising in this way commute. This
potentially leaves many diagrams in the category of interest that do not
commute. Furthermore, the interpreting functor is not necessarily
faithful. The advantages of this approach are manifold. The coherent
diagrams are represented by rooted planar binary trees with levels and
formal primitive operations on these trees. There is a close connection
with permutations and linear orderings. This avoids Catalan numbers and, 
in my view, simplifies the combinatorics. Moreover, this provides a link 
between assocaitivity and transpositions.

Premonoidal category structure has an important role to play in the
exploration of non--associative particle statistics in quantum theories,
see Joyce \cite{wjI, wjII, wjIII, wjIV, wjV}. Recently an origin has
been suggested for quark state confinement utillising a premonoidal
statistic for $SU(3)$ colour, see Joyce \cite{wj1}. Related to this is the
issue of coupling (of quantum states) which has ambiguity. For example,
in the expression $(a\otimes b)\otimes (c\otimes d)$ there is an
ambiguity concerning whether $a$ and $b$ are coupled to form $a\otimes
b$, before or after coupling $c$ and $d$ to form $c\otimes d$. In a
premonoidal category the distinction is taken into account. This removes
all freedom in the projection of diagrams in the Racah--Wigner calculus,
see Joyce \cite{wj2}. The coherence groupoid represents the coupling of
the category. It is important to note that the premonoidal category may
in fact admit a monoidal structure. Nevertheless, the structure of
interest, representing non--associative particle statistics or an
unambiguous coupling scheme, is premonoidal. It may well be that there
exist very few premonoidal categories that do not admit a monoidal
structure.

Monoidal categories were explicitly defined by Benabou \cite{jb} and Mac
Lane \cite{sm2, sm}. The monoidal category structure is found in many
areas of physics.  In quantum groups and knot theory \cite{ck}, the
Racah--Wigner calculus \cite{wjpbhr, wj} and Feynman diagrams. The
notion of coherence has its origin in Mac Lane \cite{sm2} with the
modifications of Kelly \cite{gk} and in Stasheff \cite{js}. The original
work studied natural isomorphisms for associativity, a symmetric
commutativity and identity. These were extend to cover distributivity by
Kelly \cite{gk2} and Laplaza \cite{ml}.  The non--symmetric or braided
commutativity was studied in Joyal and Street \cite{ajrs}. This paper
re--examines natural associativity but without the Pentagon condition.
An alternative and entirely different account by Yanofsky \cite{ny} was
brought to my attention by Prof. Ross Street during the final stages of
preparing this paper. The Yanofsky approach is based on higher
dimensional category theory. In contrast this approach is based on
binary trees and a natural automorphism, $\qss $, accounting for the
non--commutativity of the Pentagon diagram. Further to this structure we
may incorporate to $\qss $--Square diagrams giving a pseudomonoidal
structure. One may think of the automorphism $\qss $ as a deformation of
the Pentagon diagram. The power of this approach is realised by the
binary tree representation of coherent diagrams. The natural
automorphism $\qss $ has the simple interpretation of interchanging the
level of internal nodes. This insight suggests the obvious extension to
interchange of terminal nodes called $\qss $--pseudomonoidal.

The incorporation of the notion of identity has proven to be a delicate
balance between pseudomonoidal and $\qss $--pseudomonoidal structures.
This intermediate structure is called $\qss $--braided pseudomonoidal.
One can always account for identities by imposing the Triangle diagram.
This ultimately conflicts with the motivation behind this paper.
Although the $\qss $--braided pseudomonoidal structure carries a true
$\qss$--identity structure, coherence results from a finite number of
diagrams only in the presences of a symmetric commutativity.

In section two we define the notion of binary tree required for what
follows.  Also we define the notion of coherence for premonoidal
structures. Section three defines a premonoidal and pseudomonoidal
categories, introduces $\qss $, the coherence groupoid and proves
coherence. These two sections illustrate the methodology underlying this
paper. In section four we spell out the link to permutations and linear
orderings and briefly discuss polytopes. Section five extends the notion
of a pseudomonoidal category to the stronger $\qss $--pseudomonoidal category.
This requires extra diagrams. Namely the Dodecagon diagram and two
Quaddecagon diagrams. The coherence theorem proved in this section is
the major proof of this paper. Section six defines $\qss $--braided
pseudomonoidal categories which relaxes the conditions of section five.
Only the Decagon diagram and $\qss $--Pentagon diagram are retained.
This section uncovers a braid structure for $\qss $ where it is revealed
that the Dodecagon diagram is a Yang--Baxter condition. We call $\qss $
a $\qss $--braid because it satisfies a braid coherence result but
differs from a usual braid in that no objects are interchanged.

In section seven we begin the quest for a $\qss $--monoidal structure by
adding identities to pseudomonoidal structures. The result is given the
prefix restricted. They are monoidal whenever the identity object
indexes the natural isomorphisms. Section eight incorporates a symmetric
commutativity natural isomorphism into a $\qss $--pseudomonoidal
category. This requires $\qss $ to be symmetric, the usual Hexagon
diagram, a square diagram and two decagon diagrams. The symmetric
pseudomonoidal category requires an additional two square diagrams. In
section nine we give what rightly deserves to be called a symmetric
$\qss $--monoidal category. This requires the large and small $\qss
$--Triangle diagrams. Section ten is a summary of the premonoidal
structures of this paper.

\section{Coherence}                                              %

This section is an outline of the notion of, and the approach taken to
coherence in this paper. Also much of the notation used in the following
sections is established here. For each premonoidal type structure
presented we need a category that formally encodes what diagrams should
commute. All such categories in this paper will be some groupoid of
rooted planar binary trees for which each node (of each tree) is
assigned a level. We begin by making this notion of planar binary tree
precise. Note that we reserve the term vertex for diagrams. Instead we
use the term node.

A prerigid binary tree is a quadtuple $B=(V,E,l,s)$ consisting of a set
of nodes $V$, a set of edges $E\subset V\times V$, a level function
$l:V\rightarrow \bN \cup \{ 0\} $ and a hand function $s:\{v\in
V:l(v)\neq 0\} \rightarrow \{l,r\} $ with the properties: $(V,E)$ is a
rooted binary tree; the tree grows upward; there are no levels skipped;
and all but the empty tree have a node at level $0$. The unique node at
level $0$ is the only valency two node and is called the root. The nodes
of valency one are called leaves. The nodes of valency three are called
branches (shorthand for branch point).  The branch nodes and the root
node are called internal nodes while the leaves are terminal nodes. The
terminates (or children) for an internal node are the two unique higher
level nodes to which it is attached. We define the height of $B$ to be
$\hght{B}=\maxm\ l(V)$. The height is the level of at least one leaf.
Note that $\log_{2} |V|\leq \hght{B} \leq |V|$. The hand function $s$
assigns a left hand ($l$) or right hand ($r$) side to each terminate at
any given internal node. Two prerigid binary trees are isomorphic if
there is a bijective function between their respective node sets
preserving edges and the level and hand maps.

A rigid binary tree is an isomorphism class of prerigid binary trees.
Thus every rigid binary tree is independent of an particular set of
nodes. The name rigid is justified because there is absolutely no
topological freedom in how the binary tree can be drawn in a plane
provided we stipulate that edges do not cross. Moreover, this allows us
to assign a relative position to the leaves. We number leaf positions in
order from left to right (tracing around the top boundary of the tree).
A typical example of a rigid binary is given in figure 1.
\begin{figure}[h]
\hspace*{2cm}
\epsfxsize=120pt
\epsfbox{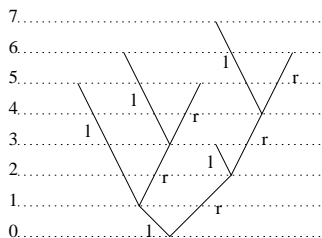}
\caption{An example of a typical rigid binary tree.}
\end{figure}
The nodes occur where lines join. We have labelled the levels and hands. There
is no need to do this and in what follows we dispense with labelling rigid
binary trees except when emphasis is required.

Let $\cC $ be a category with some premonoidal type structure. There are
many interesting and useful choices of this structure. For now we simply
accept that $\cC $ carries some such structure. We construct a free
groupoid $\cB $ over some class of rigid binary trees, together with a
length function $|\ubar |:\cB \rightarrow \bN \cup \{ 0\} $ such that
\begin{eqnarray}
\cB & = & \coprod_{n=0}^{\infty }\cB_{n}
\end{eqnarray}
where $\cB_{n}$ is the free groupoid such that $|B|=n$ for all
$B\in \cB_{n}$. The groupoid $\cB $ encodes the coherent diagrams for the
structure of $\cC $. These groupoids are not assumed to have a monoidal
structure in any sense. The objects are given by the particular class of
rigid binary trees. The arrows will be generated from a collection of
invertible primitive arrows corresponding to formal binary tree
operations. We call $\cB $ the coherence groupoid for the structure on
$\cC $.

All the diagrams we can construct in the category $\cB $ will underlie
diagrams in $\cC $ that commute. A general diagram in $\cB $ is any
closed finite directed graph. The boundary encircling each enclosed
region defines a polygonal directed graph. Every general diagram is
equivalent to the collection of polygonal directed graphs defined by its
regions. We use the word diagram in this paper to mean a polygonal
directed graph. Equivalently, a diagram is a functor into $\cB $ from the
category of two objects and two parallel arrows. The collection of all
diagrams for $\cB $ is denoted $\cD (\cB )$. Moreover,
\begin{eqnarray}
\cD (\cB ) & = & \coprod_{n=0}^{\infty }\cD (\cB_{n})
\end{eqnarray}

The relationship between $\cD (\cB )$ and the (expected) coherent diagrams in
$\cC $ is functorial. The functor of interested is the functor
\begin{eqnarray}
\can :\cB \rightarrow \coprod^{\infty }_{n=0}[\cC^{n},\cC ]~,
\end{eqnarray}
called the canonical functor, where $\can (B)\in [\cC^{|B|},\cC ]$.
Although the explicit details depend on the particular structure of
interest, we still give an outline of its construction here. It is
defined inductively according to the length of a binary tree $B\in \cB
$. The length of the binary tree $n$ determines that $\can (B)$ is an
object in $[\cC^{n},\cC]$. We call the branches of $B$ whose levels
attain the greatest level among all branch levels couples for $B$. We
will be interested in classes of rigid trees where every branch has a
distinct level and hence every tree has a unique couple. Locate the
position of the left hand most couple in $B$. The terminates are leaves
and have positions $i$ and $i+1$ for some $i\in \bN $. Removing all
three nodes and readjusting the levels we obtain a rigid binary tree
$\overline{B}$, satisfying $|\overline{B}|=|B|$ or
$|\overline{B}|=|B|-1$. In the latter case, $\can $ must satisfy the
constraint $\can (B)=\can (\overline{B})$. In the former case we define
$\can (B)$ (inductively) by
\begin{eqnarray}
\can (B)(c_{1},...,c_{n})=\can (\overline{B})(c_{1},...,c_{i-1},c_{i}\otimes
c_{i+1},c_{i+2},...,c_{n})~.
\end{eqnarray}
where $(c_{1},...,c_{n})$ is an object or arrow from $\cC^{n}$. The arrows of
$\cB $ are mapped to iterates of the natural isomorphisms of the premonoidal
structure. Thus $\can $ extends to a mapping
\begin{eqnarray}
\can :\coprod_{n=0}^{\infty }\left( \cD (\cB_{n})\times
  \cC^{n}\right) \rightarrow \cD (\cC )
\end{eqnarray}
where $\cD (\cC )$ is the set of all closed finite directed graphs in
$\cC $ whose objects are words and arrows are evaluated iterates. The
construction of $\can $ is analogous to the process of diagram
projection in the Racah--Wigner category \cite{wj2}. We define our
notion of coherence as follows.
\begin{defn}
$\cC $ is $\cB $--coherent if every diagram $D\in \cD (\cB )$ gives a
commutative diagram $\can (B)(c_{1},...,c_{|B|})$ in $\cC $ for all
objects $c_{1},...,c_{|B|}$ of $\cC^{|B|}$.
\end{defn}
For monoidal categories the coherence groupoid is rooted binary trees
($\BTree $) and $\can $ is that of Mac Lane \cite{sm2}.

\section{Premonoidal Categories}                                 %

We begin by considering a natural associativity isomorphism without any
conditions such as the Pentagon diagram.
\begin{defn}
A premonoidal category is a triple $(\cC ,\otimes ,\ass )$ where $\cC $ is 
a category, $\otimes :\cC \times \cC \rightarrow \cC $ is a bifunctor
and $\ass :\otimes (\otimes \times 1)\rightarrow \otimes (1\times
\otimes )$ is a natural isomorphism for associativity.
\end{defn}
The notion of a premonoidal category does not satisfy the Pentagon
diagram. Instead we define a natural automorphism which accounts for the
difference in the two sides. This amounts to introducing a sixth side
which could be inserted anywhere. For reasons which will reveal
themselves shortly we define the natural automorphism $\qss :\otimes
(\otimes \times \otimes ) \rightarrow \otimes (\otimes \times \otimes )$
according to the hexagonal diagram of figure 2. This is given by
composing around the bottom five sides. We call this diagram the $\qss
$--Pentagon diagram. If you set $\qss =1$ you obtain the Pentagon
diagram. In the vain of quantum groups \cite{ck} one could think of this
as a deformation of the Pentagon diagram.
\begin{figure}[h]
$\begin{diagram}
\putdtriangle<-1`0`0;600>(600,0)[(\alpha \otimes \beta )\otimes (\gamma
\otimes \delta )`((\alpha \otimes \beta )\otimes \gamma )\otimes \delta `;
\ass_{\alpha \otimes \beta ,\gamma , \delta }``]
\putqtriangle<0`1`0;600>(600,-600)[\ph{((\alpha \otimes \beta )\otimes \gamma
)\otimes \delta }``(\alpha \otimes (\beta \otimes \gamma ))\otimes \delta ;`
\ass_{\alpha ,\beta ,\gamma }\otimes 1_{\delta }`]
\puthmorphism(1200,600)[\ph{(\alpha \otimes \beta )\otimes (\gamma \otimes
\delta )}`(\alpha \otimes \beta )\otimes (\gamma \otimes \delta )`\qss_{\alpha
,\beta ,\gamma ,\delta }]{1200}{1}a
\puthmorphism(1200,-600)[\ph{(\alpha \otimes (\beta \otimes \gamma )\otimes
\delta }`\alpha \otimes ((\beta \otimes \gamma )\otimes \delta )`\ass_{\alpha
,\beta \otimes \gamma ,\delta }]{1200}{1}b
\putbtriangle<0`1`0;600>(2400,0)[\ph{(\alpha \otimes \beta )\otimes (\gamma
\otimes \delta )}``\alpha \otimes (\beta \otimes (\gamma \otimes \delta ));
`\ass_{\alpha ,\beta ,\gamma \otimes \delta }`]
\putptriangle<0`0`-1;600>(2400,-600)[`\ph{\alpha \otimes (\beta \otimes
(\gamma \otimes \delta ))}`\ph{\alpha \otimes ((\beta \otimes \gamma )\otimes
\delta )};``1_{\alpha }\otimes \ass_{\beta ,\gamma ,\delta }]
\end{diagram}$
\caption{The $\qss $--Pentagon diagram}
\end{figure}
We also have a strong version of a premonoidal category. In this case we
require $\qss_{\alpha ,\beta ,\gamma ,\delta }=\qss_{\alpha \otimes
  \beta ,\gamma \otimes \delta }$ for some natural automorphism $\qss
:\otimes \rightarrow \otimes $. Equivalently we impose the condition
that the $\qss $--Pentagon diagram holds for $\qss :\otimes \rightarrow
\otimes $. In this case $\qss_{\alpha ,\beta ,\gamma ,\delta
  }=\qss_{\alpha \otimes \beta ,\gamma \otimes \delta }$.

We now turn to the definition of the coherence groupoid for a premonoidal
category. The same groupoid describes the strong situation. This groupoid
reveals a role for the natural automorphism $\qss $.
Define $\IRBTree $ to be the free groupoid whose objects are internally
resolved binary trees (denoted IRB tree for short). An IRB tree $B$ is a rigid
binary tree where the internal nodes are assigned a distinct level and the
leaves are all assigned the same maximum level $\hght{B}$. The length of $B$ is defined to be the number of leaves.
This is given by $|B|=\hght{B}+1$. The internal nodes are labelled by
$0,...,\hght{B}-1$ according to their level. Moreover, an IRB tree is uniquely
represented by its internal node levels in the following way. Begin at the
left hand most leaf. Trace around the top of the tree. As each
internal node is passed at the bottom of a valley write its level down. This
produces an ordered sequence of the internal node levels that uniquely
describes the IRB tree. An example is given in figure 3. Every permutation of
$012...(n-1)$ represents a unique IRB tree of height $n$. Hence there are $n!$
IRB trees of height $n$.
\begin{figure}[h]
\hspace*{2cm}
\epsfxsize=60pt
\epsfbox{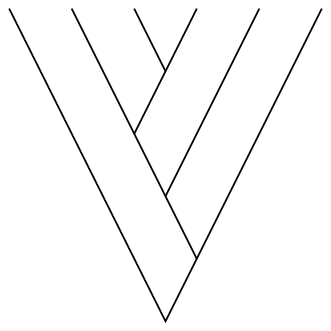}
\caption{The IRB tree described by $03421$.}
\end{figure}
The primitive arrows of $\IRBTree $ are given by formal operations on internal
nodes: One may interchange the level of a pair of adjacent branches or
reattach an edge. These formal operations are primitive
arrows corresponding to iterates of $\ass $,$\ass^{-1}$, $\qss $ and
$\qss^{-1}$; and are depicted in figure 4.
\begin{figure}[h]
\hspace*{2cm}
\epsfxsize=300pt
\epsfbox{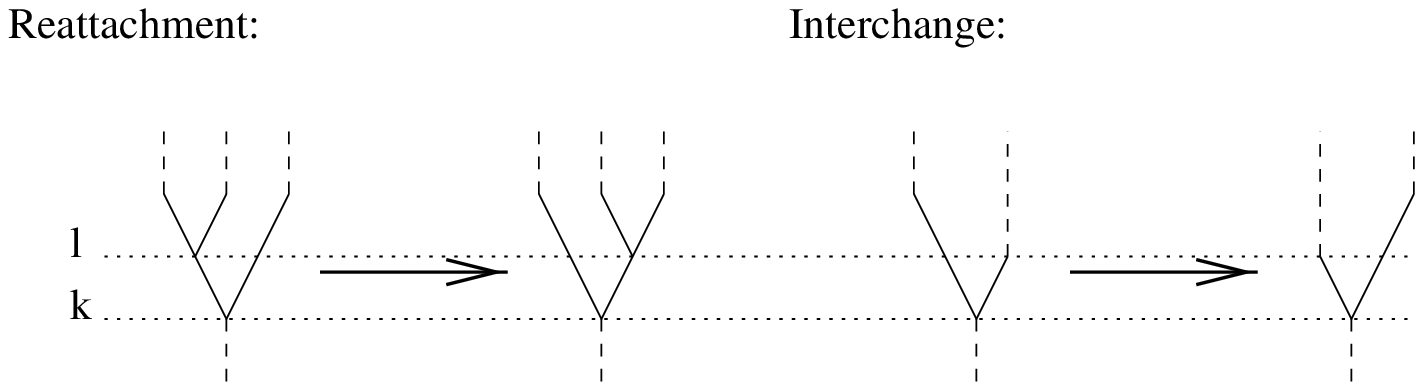}
\caption{The primitive arrows of $\IRBTree $ corresponding to iterates of
$\ass $ and $\qss $.}
\end{figure}
The dashed lines represent attachment sites to the remaining edges and
nodes of the binary tree. The node at level $k$ we call the pivot of the
arrow. Recall that the two nodes joining to $k$ from above are called
terminates.  Note that for reattachment the lowest terminate level for
$k$ must be greater than $l$. We emphasize that we do not require
$l=k+1$. If the source for reattachment has a node at level $l$ on the
left (resp. right) the arrow corresponds to an iterate of $\ass $ (resp.
$\ass^{-1}$). If the source for interchange has a node at level $l$ on
the right (resp. left) the arrow corresponds to an iterate of $\qss $
(resp. $\qss^{-1}$). The $\can $ functor is generated by $\can \ass
=\ass $, $\can \ass^{-1}=\ass^{-1}$, $\can \qss =\qss $ and $\can
\qss^{-1}=\qss^{-1}$.

We introduce some useful notation at this point. Let $B, B^{\prime }$ be
IRB trees. If there is a primitive arrow pivoting about $k$ between $B$
and $B^{\prime }$ we denote the image under $\can $ by
\begin{eqnarray}
\It^{B}_{k}(\iss ):\can (B)\rightarrow \can (B^{\prime })~,
\end{eqnarray}
where $\iss $ is any natural isomorphism of the premonoidal structure.
We usually dispense with writing the superscript $B$. For example
$\It_{2}(\ass ) =1\otimes (\ass \otimes 1)$ for the IRB tree of figure
3.
\begin{prop}$\ph{nothing}$
\begin{enumerate}
\item{Given two IRB trees of height $n$ then there is a finite
    sequence of primitive arrows transforming one into the other.}
\item{Every IRB tree of height $n$ is the source of no more than $n-1$
    distinct primitive arrows.}
\item{There are $n!$ distinct IRB trees of height $n$.}
\end{enumerate}
\end{prop}
\pf\ (i) We prove by induction on $n$ that every IRB tree
of height $n$ may be brought into the form $012...(n-1)$. Consider a
binary tree $a_{0}a_{1}...a_{n}$ of height $n+1$. If $a_{0}\neq n$ then 
we may apply the induction hypothesis to $a_{0}...a_{i-1}a_{i+1}...a_{n}$
where $a_{i}=n$ is omitted. In particular there is a sequence of primitive
arrows transforming the IRB tree to $0...(i-1)(i+1)...(n-1)$. Hence
$a_{0}...a_{n}$ can be transformed to $0...(i-1)n(i+1)...(n-1)$. Again
the induction hypothesis may be applied to $1...(i-1)n(i+1)...(n-1)$ to bring
it into the form $1...n$. Hence $a_{0}...a_{n}$ may be brought into the form
$0...n$. If $a_{0}=n$ then by the induction hypothesis we can arrange the last
$n$ terms of $a_{0}...a_{n}$ as we wish and hence bring it into the form
$n(n-1)...10$. The primitive arrow given by the transposition $(12)$
transforms the tree to $(n-1)n(n-2)...10$. The first term is not $n$ so by
the first case it may be brought into the form $01...n$.\\
(ii) We prove by induction on the height $n$. Suppose the hypothesis holds for
$a_{0}...a_{n-1}$. Consider $a_{0}...a_{n}$. Ignoring $a_{i}=n$ we can
apply a maximum of $n$ distinct primitive operations by the induction
hypothesis. Any additional operations on $a_{0}...a_{n}$ involve a
transposition moving $n$. There is at most only one possible such primitive
operation. Hence there are at most $n$ distinct primitive arrows with source
$a_{0}...a_{n}$.\\
(iii) This has already been noted.

Each IRB tree $B$ gives a functor $\can (B):\cC^{|B|}\rightarrow \cC $
given by bracketing according to the binary tree. Note that different
IRB trees may map to the same objects and arrows. For example the trees
$201$ and $102$ both correspond to the functor $(\ubar \otimes \ubar
)\otimes (\ubar \otimes \ubar )$. The distinction between the two trees
corresponds to a formal weight on the brackets. That is, $201$
corresponds to $(\ubar \otimes \ubar )_{2}\otimes (\ubar \otimes \ubar
)_{1}$ and $102$ corresponds to $(\ubar \otimes \ubar )_{1}\otimes
(\ubar \otimes \ubar )_{2}$. The primitive arrows of $\IRBTree $ map to
an iterate of one of the natural isomorphisms $\ass ,\ass^{-1},\qss
,\qss^{-1}$.  The $\qss $--Pentagon diagram has the underlying $\IRBTree
$ diagram structure given by figure 5. The notation $(ij)$ means the
natural isomorphism corresponding to the transposition swaping the $i$th 
and $j$th positions in the linear ordering.
\begin{figure}[h]
\hspace*{2cm}
\epsfxsize=200pt
\epsfbox{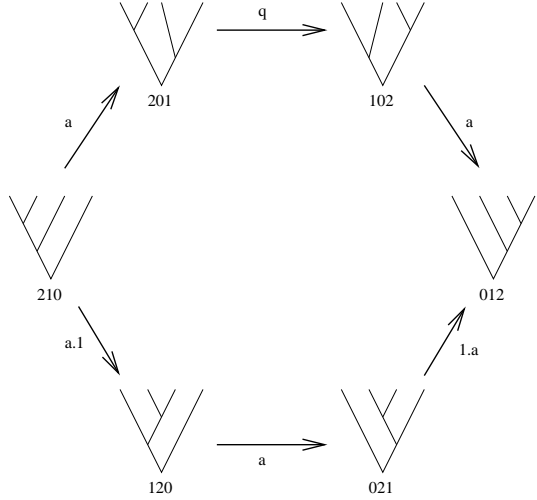}
\caption{The length four diagram in $\IRBTree $ underlying the
$\qss $--Pentagon diagram.}
\end{figure}

\begin{defn}
a pseudomonoidal category is a premonoidal category $(\cC ,\otimes ,\ass 
)$ satisfying the two $\qss $--Square diagrams of figure 6.
\end{defn}
\begin{figure}[h]
$\begin{diagram}
\putsquare<1`1`1`1;1600`600>(800,0)[((\alpha \otimes \beta )\otimes
\gamma )\otimes (\delta \otimes \epsilon )`((\alpha \otimes \beta )\otimes
\gamma )\otimes (\delta \otimes \epsilon )`(\alpha \otimes (\beta \otimes
\gamma ))\otimes (\delta \otimes \epsilon )`(\alpha \otimes (\beta \otimes
\gamma ))\otimes (\delta \otimes \epsilon );\qss_{\alpha \otimes \beta
  ,\gamma ,\delta ,\epsilon }`\ass_{\alpha ,\beta ,\gamma }\otimes
1_{\delta \otimes \epsilon }`\ass_{\alpha ,\beta ,\gamma }\otimes
1_{\delta \otimes \epsilon }`\qss_{\alpha ,\beta \otimes \gamma ,\delta
  ,\epsilon }]
\end{diagram}$\\ \\ \\
$\begin{diagram}
\putsquare<1`1`1`1;1600`600>(800,0)[(\alpha \otimes \beta )\otimes
((\gamma \otimes \delta )\otimes \epsilon )`(\alpha \otimes \beta )\otimes
((\gamma \otimes \delta )\otimes \epsilon )`(\alpha \otimes \beta )\otimes
(\gamma \otimes (\delta \otimes \epsilon ))`(\alpha \otimes \beta )\otimes
(\gamma \otimes (\delta \otimes \epsilon ));\qss_{\alpha ,\beta
  ,\gamma \otimes \delta ,\epsilon }`1_{\alpha \otimes \beta }\otimes
\ass_{\gamma ,\delta ,\epsilon }`1_{\alpha \otimes \beta }\otimes
\ass_{\gamma ,\delta ,\epsilon }`\qss_{\alpha ,\beta ,\gamma ,\delta
  \otimes \epsilon }]
\end{diagram}$
\caption{The $\qss $--Square diagrams.}
\end{figure}
If one substitutes for deformativity using the $\qss $--Pentagon in the
above $\qss $--Square diagrams then one sees that these diagrams are
actually dodecagons. In this situation the primitive arrows of figure 4
relax the adjacent level requirement. Thus the two $\qss $--Square
diagrams correspond to the $\IRBTree $ digrams in figure 7.
\begin{figure}[h]
\hspace*{2cm}
\epsfxsize=320pt
\epsfbox{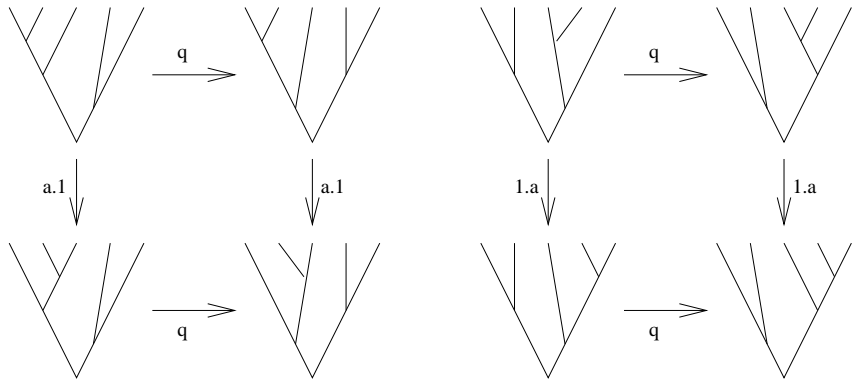}
\caption{The length five diagrams in $\IRBTree $ underlying the two
$\qss $--Square diagrams.}
\end{figure}

We turn now to the coherence of these categories. For a Premonoidal
category every $\qss $ is defined by a sequence of five reattachment
arrows under $\can $. The $\qss $ automorphisms may be factored out.
Thus all diagrams induced under $\can $ commute as a result of natural
diagrams.
\begin{thm}
  Every pseudomonoidal (resp. strong pseudomonoidal) category $(\cC
  ,\otimes ,\ass )$ is pseudomonoidal coherent (resp. strong
  pseudomonoidal) coherent if and only if the $\qss $--Pentagon diagram
  and two $\qss $--Square diagrams hold.
\end{thm}
\pf\ Let $D$ be a diagram in $\IRBTree $. Note that each vertex has
the same length. Let this length be $n$. The rank of a diagram is defined to
be the length of any one of its vertices. We prove coherence by
induction on diagram rank. The result can be verified explicitly for the
ranks $1,2,3,4$. The rank four case is given in Appendix A. Suppose the
result holds for diagrams of rank $n+1$. Let $a_{0},...,a_{r}$ be the
vertices for some diagram $D$ of rank $n+2$ given by reading around the
outside. We identify $a_{r+1}$ with $a_{0}$. If $n$ always occurs in the
first position of vertex $a_{i}$ then $D$ commutes by naturality
and the induction hypothesis. Now suppose that $n$ does not occur in the
first position of some vertices of $D$. We divide $D$ into alternate
maximal sections where $n$ is in the first position alternating with $n$
is never in the first position. We replace, using the Pentagon diagram,
all arrows raising/lowering the first position branch to/from level
$n$. Now we can and do assume that every arrow moving $n$ into or out of
the first position corresponds to an iterate of $\ass^{-1}$ or $\ass $
respectively.

Let a typical maximal section with $n$ in the first position be $a_{i}
\rightarrow \cdots \rightarrow a_{j}$. Let the arrows $a_{i-1}\rightarrow
a_{i}$ and $a_{j}\rightarrow a_{j+1}$ be $\It_{k}(\ass^{-1})$ and
$\It_{l}(\ass )$ respectively. We show how to replace the sequence $a_{i-1}
\rightarrow \cdots \rightarrow a_{j+1}$ with an alternative sequence where the
first position is never $n$. Moreover, we assume that $k<l$. The modification
to the other case is obvious. The construction is depicted in figure 8.
\begin{figure}[h]
\hspace*{2cm}
\epsfxsize=300pt
\epsfbox{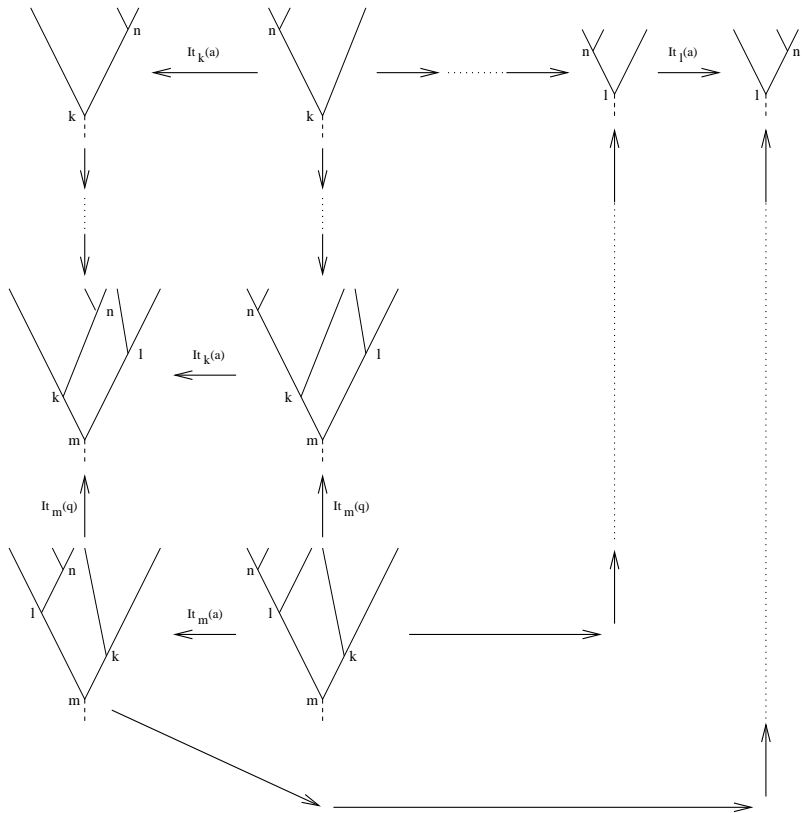}
\caption{Removal of maximal sequence with $n$ in the first position.}
\end{figure}
The sequence of arrows along the top is $a_{i-1}\rightarrow \cdots
\rightarrow a_{j+1}$. The vertical sides of the top left region are
identical and keep the position of $n$ fixed. The bottom arrow of this
region moves $n$ into the first position. This region commutes by the
induction hypothesis, $\qss $--Square diagrams and naturality. The next
region is a natural square that interchanges the levels $k$ and $l$. The
bottom and right edges of the centre region is a sequence of arrows
keeping $n$ and $l$ fixed. The region enclosed always keeps $n$ fixed
and so by hypothesis commutes. Finally the sequence around the bottom
going up the right hand side is a sequence of arrows with $n$ and $l$
fixed enclosing the last region. This region commutes by naturality and
the induction hypothesis.

If $n$ is not in a fixed position for every vertex of $D$, then we apply
the above argument to exclude $n$ from position two (as well as position
one).  We repeat this argument inductively for position three, and so
on, only stopping when $n$ is in a fixed position for all vertices of
$D$. The diagram $D$ now commutes by the induction hypothesis. This
completes the proof.

\section{Natural Associativity, Permutations and Linear Orderings}%

Any natural isomorphism for associativity has a close relationship with
the Symmetric groups. The $\qss $ automorphism that arises accounts for
the degenerate nature of the functor $\otimes (\otimes \times \otimes
):\cC^{4}\rightarrow \cC $. The Pentagon diagram reflects this
degeneracy. The distinct is made in the $\qss $--Pentagon diagram. Even
if the Pentagon diagram holds the distinction can always be made at the
formal level. In a nutshell a premonoidal structure allows one to
utilise the symmetric groups and avoid the combinatorics of Catalan
numbers. We spell out the precise connection here. A similar
relationship, this time between the associahedron and the permutohedron,
has been given in Tonks \cite{at}.

Let $\cLO_{n}$ be the groupoid of linear orderings of length $n$. The
objects are linear orderings of $0,1,...,n-1$ and the arrows are
permutations. The groupoid of all linear orderings is given by $\cLO =
\coprod_{n=0}^{\infty } \cLO_{n}$. Let $F:\IRBTree \rightarrow \cLO $ be
the functor outlined in the previous section. Each IRB tree is mapped to
the sequence of its internal levels and the arrows are mapped to
transpositions. Premonoidal and pseudomonoidal coherence (Theorem 1)
implies that $F$ is bijective. Moreover, we can extend the canonical
functor to $\cLO $.  This functor is
\begin{eqnarray}
\sym :\cLO \rightarrow \coprod_{n=0}^{\infty }[\cC^{n},\cC ]~,
\end{eqnarray}
where $\sym \circ F=\can $.

Briefly we consider the construction of $\qss $--associahedra. The polytope
for words of length five is given in figure 53. This planar diagram
folds into the partially-formed truncated octahedron of figure 54. Some
faces are missing or halved of this shape. Also there
are four vertices that are the source of only two primitive arrows instead of
three. This prevents the construction of a polytope. The solution to this
dilemma is to use the permutation structure. We summarise in the first few
polytopes (or permutahedra) in the following table.
\begin{table}[h]
\hspace*{2cm}
\begin{tabular}{|c|c|c|c|c|}
\hline
$n$ & source & $\sym (12)$ & $\sym (23)$ & polytope \\
\hline
$2$ & $0$ & & & point \\
\hline
$3$ & $01$ & $\ass^{-1}$ & & line segment \\
 & $10$ & $\ass $ & & \\
\hline
$4$ & $012$ & $\ass^{-1}$ & $\ass^{-1}$ & hexagon \\
& $102$ & $\ass $ & $\ass^{-1}(1.\ass^{-1})\ass $ & \\
& $021$ & $\ass (\ass^{-1}.1)\ass^{-1}$ & $\ass $ & \\
& $201$ & $\ass (\ass .1)\ass^{-1}$ & $\ass^{-1}$ & \\
& $210$ & $\ass$ & $\ass $ & \\
& $120$ & $\ass^{-1}$ & $\ass^{-1}(1.\ass )\ass $ & \\
\hline
$5$ & $0123$ & $\ass^{-1}$ & $1.\ass^{-1}$ & truncated octahedron \\
& \vdots &\vdots &\vdots & \\
\hline
\end{tabular}
\caption{The generators of $\cLO_{n}$ where $2\leq n\leq 5$.}
\end{table}
Note that because of the $\qss $--Pentagon diagram every permutation under
$\sym $ can be described by a sequence of iterates of $\ass $.
\begin{figure}[h]
\hspace*{2cm}
\epsfxsize=200pt
\epsfbox{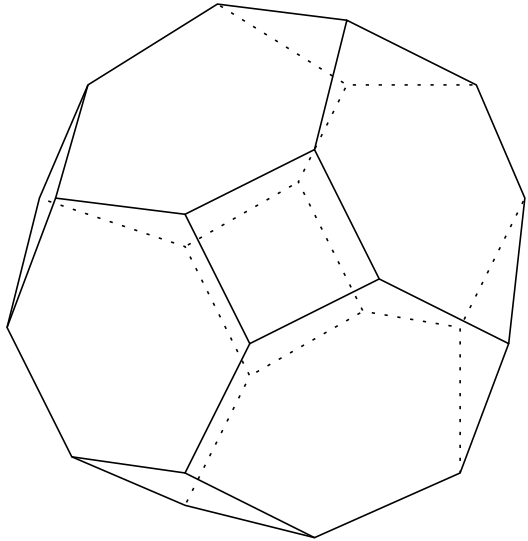}
\caption{A truncated octahedron.}
\end{figure}

\section{$\qss $--Pseudomonoidal Categories}                     %

The $\qss $ natural automorphism for pseudomonoidal categories was found to
behave (at the formal level at least) as an interchange of the internal
node level. In this section we extend this idea to
include the interchange of leaf levels as well.
\begin{defn}
  A $\qss $--pseudomonoidal category is a quadtruple $(\cC ,\otimes
  ,\ass ,\qss )$ where $\cC $ is a category, $\otimes :\cC \times \cC
  \rightarrow \cC $ is a bifunctor, $\ass :\otimes (\otimes \times
  1)\rightarrow \otimes (1\times \otimes )$ and $\qss :\otimes
  \rightarrow \otimes $ are natural isomorphisms satisfying $(\cC
  ,\otimes ,\ass ,\qss )$ is strong pseudomonoidal and the Dodecagon
  diagram (figure 10) and the two Quaddecagon diagrams (figures 11 and 12)
  hold.
\end{defn}
Note that we sometimes use a dot as an abbreviation of $\otimes $ (as in
figures 11 and 12) and often dispense with subscripts on the natural
isomorphisms.
\begin{figure}[h]
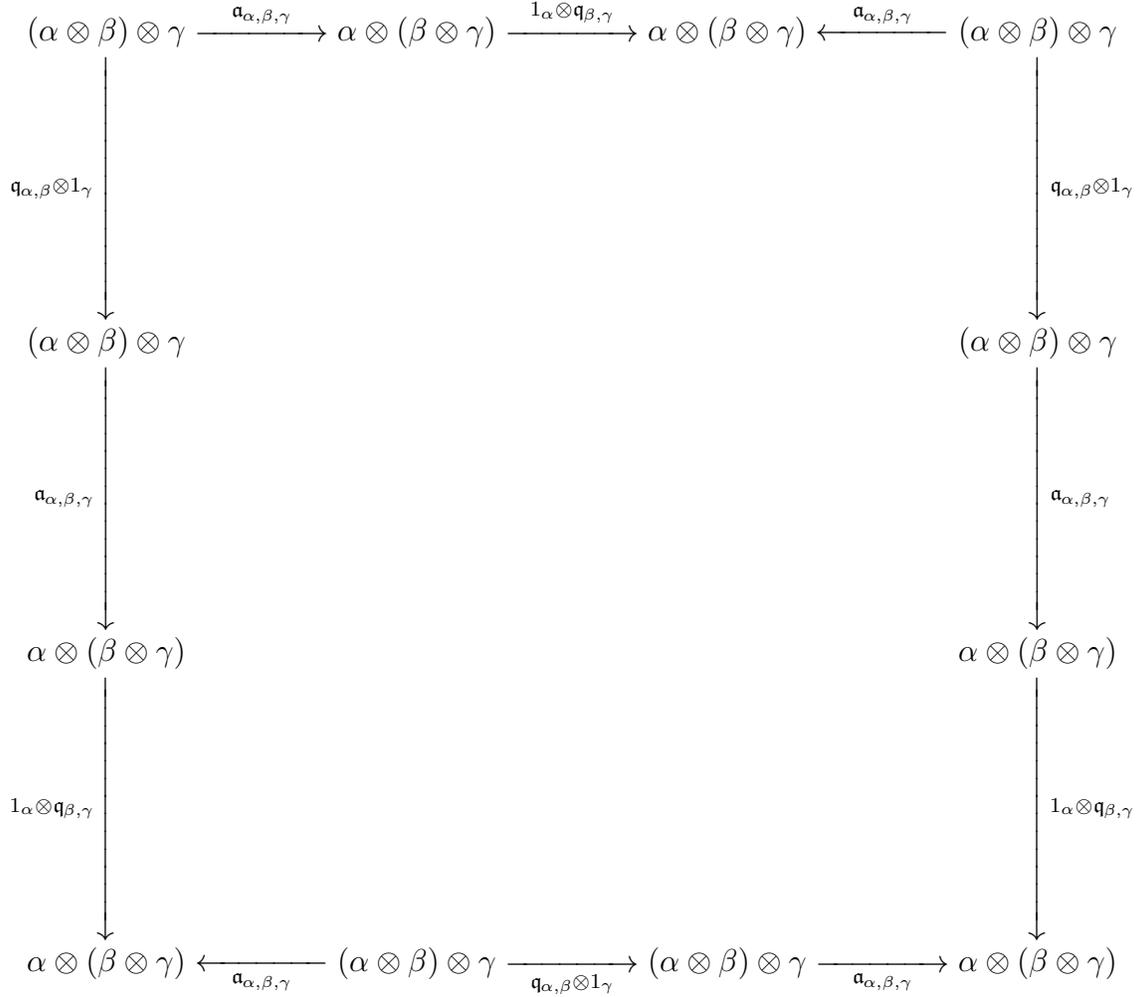

$\begin{diagram}
\puthmorphism(500,0)[(\alpha \otimes \beta )\otimes \gamma `\alpha \otimes
(\beta \otimes \gamma )`\ass_{\alpha ,\beta ,\gamma }]{1000}{1}a
\puthmorphism(1500,0)[\ph{\alpha \otimes (\beta \otimes \gamma )}`\ph{\alpha
\otimes (\beta \otimes \gamma )}`1_{\alpha }\otimes \qss_{\beta ,
\gamma }]{1000}{1}a
\puthmorphism(2500,0)[\alpha \otimes (\beta \otimes \gamma )`(\alpha \otimes
\beta )\otimes \gamma `\ass_{\alpha ,\beta ,\gamma }]{1000}{-1}a
\putvmorphism(500,0)[\ph{(\alpha \otimes \beta )\otimes \gamma }`\ph{(\alpha
\otimes \beta )\otimes \gamma }`\qss_{\alpha ,\beta }\otimes 1_{\gamma }]
{1000}{1}l
\putvmorphism(500,-1000)[(\alpha \otimes \beta )\otimes \gamma `\alpha \otimes
(\beta \otimes \gamma )`\ass_{\alpha ,\beta ,\gamma }]{1000}{1}l
\putvmorphism(500,-2000)[\ph{\alpha \otimes (\beta \otimes \gamma )}`
\ph{\alpha \otimes (\beta \otimes \gamma )}`1_{\alpha }\otimes \qss_{\beta 
,\gamma }]{1000}{1}l
\putvmorphism(3500,0)[\ph{(\alpha \otimes \beta )\otimes \gamma }`\ph{(\alpha
\otimes \beta )\otimes \gamma }`\qss_{\alpha ,\beta }\otimes 1_{\gamma }]
{1000}{1}r
\putvmorphism(3500,-1000)[(\alpha \otimes \beta )\otimes \gamma `\alpha \otimes
(\beta \otimes \gamma )`\ass_{\alpha ,\beta ,\gamma }]{1000}{1}r
\putvmorphism(3500,-2000)[\ph{\alpha \otimes (\beta \otimes \gamma )}`
\ph{\alpha \otimes (\beta \otimes \gamma )}`1_{\alpha }\otimes \qss_{\beta 
,\gamma }]{1000}{1}r
\puthmorphism(500,-3000)[\alpha \otimes (\beta \otimes \gamma )`(\alpha \otimes
\beta )\otimes \gamma `\ass_{\alpha ,\beta ,\gamma }]{1000}{-1}b
\puthmorphism(1500,-3000)[\ph{(\alpha \otimes \beta )\otimes \gamma }`
\ph{(\alpha \otimes \beta )\otimes \gamma }`\qss_{\alpha ,\beta }\otimes
1_{\gamma }]{1000}{1}b
\puthmorphism(2500,-3000)[(\alpha \otimes \beta )\otimes \gamma `\alpha
\otimes (\beta \otimes \gamma )`\ass_{\alpha ,\beta ,\gamma }]{1000}{1}b
\end{diagram}$
\caption{The Dodecagon diagram.}
\end{figure}
\begin{figure}[h]
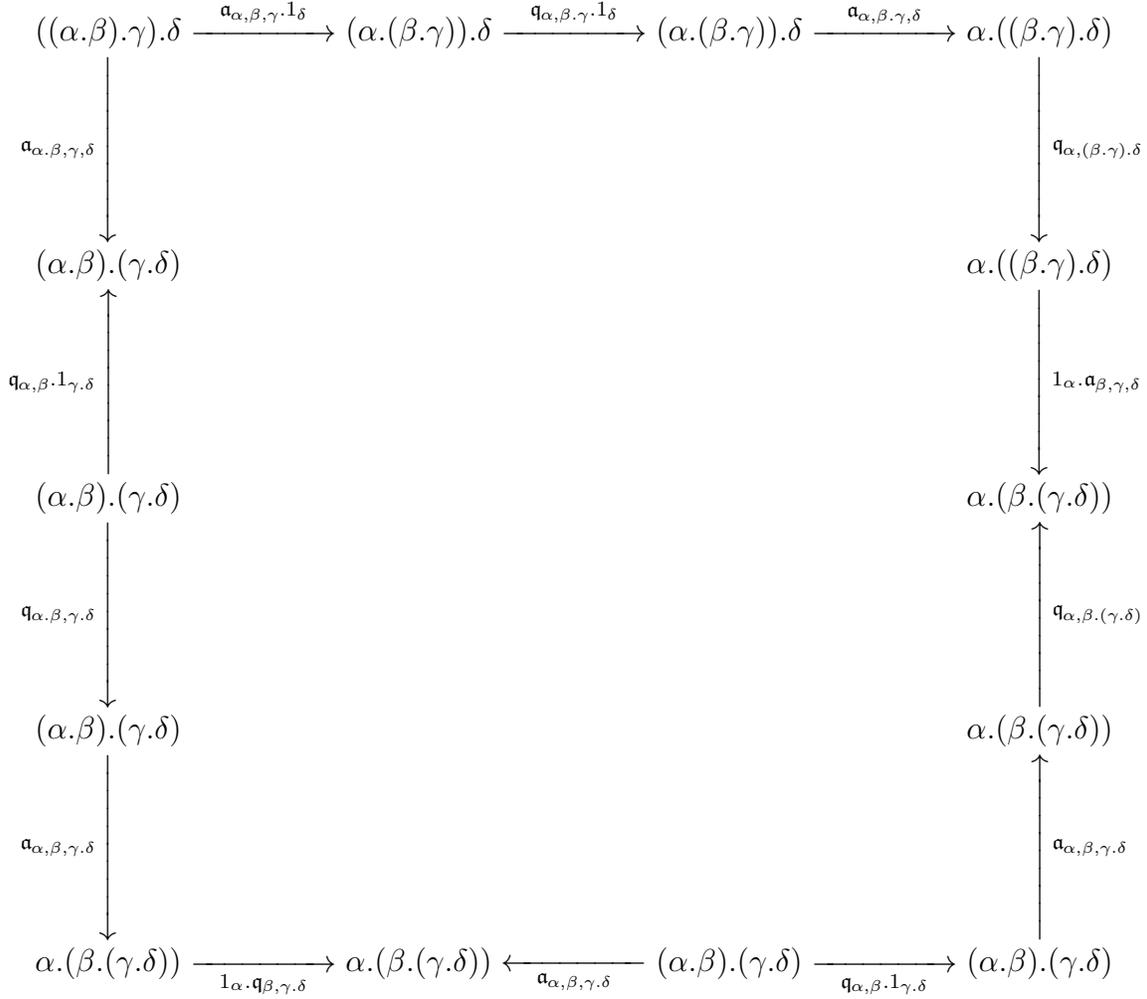

$\begin{diagram}
\puthmorphism(500,0)[((\alpha .\beta ).\gamma ).\delta `(\alpha .(\beta .\gamma
)).\delta `\ass_{\alpha ,\beta ,\gamma }.1_{\delta }]{1000}{1}a
\puthmorphism(1500,0)[\ph{(\alpha .(\beta .\gamma )).\delta }`(\alpha .(\beta
.\gamma )).\delta `\qss_{\alpha ,\beta .\gamma }.1_{\delta }]{1000}{1}a
\puthmorphism(2500,0)[\ph{(\alpha .(\beta .\gamma )).\delta }`\alpha .((\beta
.\gamma ).\delta )`\ass_{\alpha ,\beta .\gamma ,\delta }]{1000}{1}a
\putvmorphism(500,0)[\ph{((\alpha .\beta ).\gamma ).\delta }`(\alpha .\beta )
.(\gamma .\delta )`\ass_{\alpha .\beta ,\gamma ,\delta }]{750}{1}l
\putvmorphism(500,-750)[\ph{(\alpha .\beta ).(\gamma .\delta )}`(\alpha .\beta
).(\gamma .\delta )`\qss_{\alpha ,\beta }.1_{\gamma .\delta }]{750}{-1}l
\putvmorphism(500,-1500)[\ph{(\alpha .\beta ).(\gamma .\delta )}`(\alpha
.\beta ).(\gamma .\delta )`\qss_{\alpha .\beta ,\gamma .\delta }]{750}{1}l
\putvmorphism(500,-2250)[\ph{(\alpha .\beta ).(\gamma .\delta )}`\alpha
.(\beta .(\gamma .\delta ))`\ass_{\alpha ,\beta ,\gamma .\delta }]{750}{1}l
\puthmorphism(500,-3000)[\ph{\alpha .(\beta .(\gamma .\delta ))}`\alpha
.(\beta .(\gamma .\delta ))`1_{\alpha }.\qss_{\beta ,\gamma .\delta }]{1000}
{1}b
\puthmorphism(1500,-3000)[\ph{\alpha .(\beta .(\gamma .\delta ))}`(\alpha
.\beta ).(\gamma .\delta )`\ass_{\alpha ,\beta ,\gamma .\delta }]{1000}{-1}b
\puthmorphism(2500,-3000)[\ph{(\alpha .\beta ).(\gamma .\delta )}`(\alpha
.\beta ).(\gamma .\delta )`\qss_{\alpha ,\beta }.1_{\gamma .\delta }]{1000}{1}b
\putvmorphism(3500,0)[\ph{\alpha .((\beta .\gamma ).\delta )}`\alpha .((\beta
.\gamma ).\delta )`\qss_{\alpha ,(\beta .\gamma ).\delta }]{750}{1}r
\putvmorphism(3500,-750)[\ph{\alpha .((\beta .\gamma ).\delta )}`\alpha
.(\beta .(\gamma .\delta ))`1_{\alpha }.\ass_{\beta ,\gamma ,\delta }]{750}{1}r
\putvmorphism(3500,-1500)[\ph{\alpha .(\beta .(\gamma .\delta ))}`\alpha
.(\beta .(\gamma .\delta ))`\qss_{\alpha ,\beta .(\gamma .\delta )}]{750}{-1}r
\putvmorphism(3500,-2250)[\ph{\alpha .(\beta .(\gamma .\delta ))}`\ph{(\alpha
.\beta ).(\gamma .\delta )}`\ass_{\alpha ,\beta ,\gamma .\delta }]{750}{-1}r
\end{diagram}$
\caption{The first Quaddecagon diagram.}
\end{figure}
\begin{figure}[h]
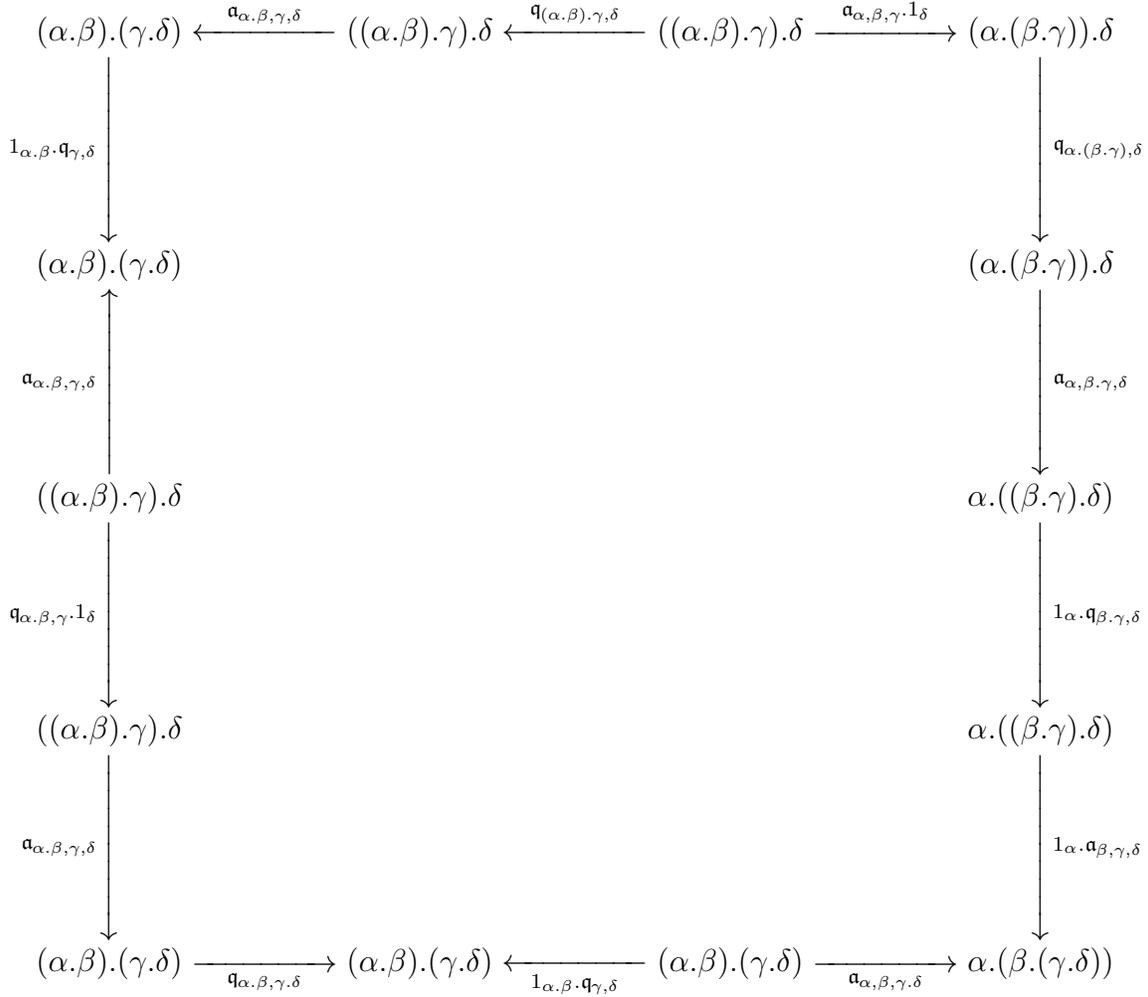

$\begin{diagram}
\puthmorphism(500,0)[(\alpha .\beta ).(\gamma .\delta )`((\alpha .\beta )
.\gamma ).\delta `\ass_{\alpha .\beta ,\gamma ,\delta }]{1000}{-1}a
\puthmorphism(1500,0)[\ph{((\alpha .\beta ).\gamma ).\delta }`((\alpha .\beta
).\gamma ).\delta `\qss_{(\alpha .\beta ).\gamma ,\delta }]{1000}{-1}a
\puthmorphism(2500,0)[\ph{((\alpha .\beta ).\gamma ).\delta }`(\alpha .(\beta
.\gamma )).\delta `\ass_{\alpha ,\beta ,\gamma }.1_{\delta }]{1000}{1}a
\putvmorphism(500,0)[\ph{(\alpha .\beta ).(\gamma .\delta )}`(\alpha .\beta )
.(\gamma .\delta )`1_{\alpha .\beta }.\qss_{\gamma ,\delta }]{750}{1}l
\putvmorphism(500,-750)[\ph{(\alpha .\beta ).(\gamma .\delta )}`((\alpha .\beta
).\gamma ).\delta `\ass_{\alpha .\beta ,\gamma ,\delta }]{750}{-1}l
\putvmorphism(500,-1500)[\ph{((\alpha .\beta ).\gamma ).\delta }`((\alpha
.\beta ).\gamma ).\delta `\qss_{\alpha .\beta ,\gamma }.1_{\delta }]{750}{1}l
\putvmorphism(500,-2250)[\ph{((\alpha .\beta ).\gamma ).\delta }`(\alpha
.\beta ).(\gamma .\delta )`\ass_{\alpha .\beta ,\gamma ,\delta }]{750}{1}l
\puthmorphism(500,-3000)[\ph{(\alpha .\beta ).(\gamma .\delta )}`(\alpha
.\beta ).(\gamma .\delta )`\qss_{\alpha .\beta ,\gamma .\delta }]{1000}{1}b
\puthmorphism(1500,-3000)[\ph{(\alpha .\beta ).(\gamma .\delta )}`(\alpha
.\beta ).(\gamma .\delta )`1_{\alpha .\beta }.\qss_{\gamma ,\delta }]{1000}
{-1}b
\puthmorphism(2500,-3000)[\ph{(\alpha .\beta ).(\gamma .\delta )}`\alpha
.(\beta .(\gamma .\delta ))`\ass_{\alpha ,\beta ,\gamma .\delta }]{1000}{1}b
\putvmorphism(3500,0)[\ph{(\alpha .(\beta .\gamma )).\delta }`(\alpha .(\beta
.\gamma )).\delta `\qss_{\alpha .(\beta .\gamma ),\delta }]{750}{1}r
\putvmorphism(3500,-750)[\ph{(\alpha .(\beta .\gamma )).\delta }`\alpha
.((\beta .\gamma ).\delta )`\ass_{\alpha ,\beta .\gamma ,\delta }]{750}{1}r
\putvmorphism(3500,-1500)[\ph{\alpha .((\beta .\gamma ).\delta )}`\alpha
.((\beta .\gamma ).\delta )`1_{\alpha }.\qss_{\beta .\gamma ,\delta }]{750}{1}r
\putvmorphism(3500,-2250)[\ph{\alpha .((\beta .\gamma ).\delta )}`\ph{\alpha
.(\beta .(\gamma .\delta ))}`1_{\alpha }.\ass_{\beta ,\gamma ,\delta }]{750}
{1}r
\end{diagram}$
\caption{The second Quaddecagon diagram.}
\end{figure}

In order to understand the Dodecagon and Quaddecagon diagrams we need to
understand the underlying combinatorics. Define a resolved binary tree
or RB tree to be a rigid binary tree where every node is assigned a
distinct level. We represent an RB tree by a finite sequence of levels
as follows. Tracing around the top of the tree beginning at the left
hand side leaf, we generate a sequence of all the node levels,
$a_{1},a_{2},...,a_{2n-1}$, where $n$ is the length of the tree. We also
write this sequence in the exploded form
\begin{eqnarray}
\begin{array}{ccccccccc}a_{1} & & a_{3} & & a_{5} & \cdots & a_{2n-3}
& & a_{2n-1} \\ & a_{2} & & a_{4} & & \cdots & & a_{2n-2} & \end{array}~.
\end{eqnarray}
The bottom row contains the internal node levels and the top the leaf levels.
The levels satisfy $a_{2i}<a_{2i-1},a_{2i+1}$ for all $i=1,2,...,n-1$.
Moreover, any such format of the numbers $0,1,...,2n-1$ obeying this condition
uniquely determines an RB tree. We let $\RBTree $ denote the free groupoid
of RB trees. The objects are RB trees. The arrows are generated from the
primitive operations for reattachment of adjacent edges (corresponding to
associativity) and interchange of level (corresponding to $\qss $). These
operations are the obvious extensions of the primitive operations of
$\IRBTree $ (figure 4). The degree of connectedness and size of $\RBTree $ is
given by Proposition 2.
\begin{prop}$\ph{nothing}$
\begin{enumerate}
\item{Given two RB trees of length $n$ then there is a finite sequence
of primitive arrows transforming one into the other.}
\item{Every RB tree of length $n$ is the source of no more than
$2(n-1)$ distinct primitive arrows.}
\item{Let $T(m)$ be the number of RB trees of length $m$. $T(n)$
satisfies the recursive formula
\begin{eqnarray}
T(n)=\sum_{m=1}^{n-1}\Comb{2n-2}{2m-1}T(m)T(n-m)~, \ph{space}T(1)=1 ~.
\end{eqnarray}}
\end{enumerate}
\end{prop}
\pf\ (i) We use the principle of strong induction. Clearly the result holds for
$n=2$. Suppose it holds for $n$ and consider a RB tree $B$ of length
$n+1$. We will show that it can be brought into the standard form
$\COLbin{1\ph{0}3\ph{2}5\ph{4}\cdots 2n-2\ph{2n-3}2n-1}{\ph{1}0\ph{3}2\ph{5}4
\cdots \ph{2n-2}2n-3\ph{2n-2}}$. It follows from this, because the
arrows are invertible, that any two trees that can be brought into standard
form can also be transformed one to the other. If $B$ begins $\COLbin{1\ph{0}
\cdots }{\ph{1}0\cdots }$ then the remaining portion may be transformed into
standard form by the induction hypothesis. If $B$ ends $\COLbin{\cdots
\ph{0}1}{\cdots 0\ph{1}}$ then $B$ may be brought into the form $\COLbin{3
\ph{2}5\ph{4}\cdots 2n-2\ph{2n-3}2n-1\ph{0}1}{\ph{3}2\ph{5}4\cdots \ph{2n-2}
2n-3\ph{2n-1}0\ph{1}}$. Interchanging $1$ and $2$ then $B$ falls in the next
case to be considered.

The only other case that can occur is when $0$ and $1$
are both in the bottom row of $B$. We suppose that $1$ is to the right of the
$0$. If not then one simply interchanges them and continues. Let
$m$ be the smallest level to the left of $0$. Hence $1,...,m-1$ are to the
right of $0$. If $m<2n+1$ then $m+1$ occurs in $B$. If it occurs to the right
of $0$ then by the induction hypothesis applied to $1,...,m-1$ we bring
$B$ into the form $\cdots m\cdots 0(m+1)1\cdots $. To this we apply the
sequence of primitive operations: interchange $0$ and $1$, then
interchange $m$ and $m+1$, and finally interchange $0$ and $1$. Thus $B$ is
in the form $\cdots (m+1)\cdots 0m1\cdots $. Hence the smallest level to the
left of $0$ is increased to $m+1$. In the other situation $m+1$ is to the
left of $0$ and on the bottom row. Applying the induction hypothesis to
$1,...,m-1$ we bring $B$ into the form $\cdots m\cdots 0\cdots (m-1)\cdots
(m-2) \cdots \cdots 2\cdots 1\cdots $ where $1,...,m-1$ are internal (and on
the bottom row in the exploded form). Next we apply the sequence of primitive
operations: Interchange $0$ with $1$, then interchange $1$ with $2$,
continuing this process until the interchange $m-2$ with $m-1$. Now $B$ is in
the form $\cdots m\cdots (m-1)\cdots \cdots 2\cdots 1\cdots 0\cdots $ and
$0,...,m$ are all internal. Next we reverse the sequence of primitive
operations: Interchange $m$ with $m-1$, then interchange $m-1$ with $m-2$,
continuing this process until  the interchange $0$ with $1$. Now $1,...,m$ are
all to the left of $0$ and the smallest level to the left of $0$ has been
increased to $m+1$. Hence, whenever $m<2n+1$ we continue to apply the above
procedures terminating when $m$ has been increased to $2n+1$. Hence $B$ is of
the form $\COLbin{2n+1\ph{0}\cdots }{\ph{2n+1}0\cdots }$. By the induction
hypothesis applied to $1,...,2n$ we bring $B$ into the form $\COLbin{2n+1
\ph{0}2\ph{1}\cdots }{\ph{2n+1}0\ph{2}1\cdots }$. Now we apply the primitive
operations: Interchange $0$ and $1$, then interchange $2n+1$ and $2$,
then interchange $0$ and $1$ again, and finally interchange $1$ and $2$. This
brings $B$ into the form $\COLbin{1\ph{0}2n+1\ph{1}\cdots }{\ph{1}0\ph{2n+1}2
\cdots }$ which we have already shown can be brought into standard form. This
completes the proof of (i).\\
(ii) This is proved by a strong induction argument. The result
is true for $n=2$. Suppose it is true for all RB trees of length less
than or equal to $n$. Any RB tree $B$ of length $n+1$ can be divided
at $0$ into two trees of length $m$ and $n-m$ respectively. A maximum of two
primitive operations can pivot about $0$ and applying the induction hypothesis
to both trees we get an upper bound of $2(m-1)+2(n-m)+2=2n$ for the maximum
number of primitive operations on $B$.\\
(iii) An RB tree of length $n$ has $2n-1$ levels with the first level
occupied by the root. We divide the tree into left hand and
right hand components. If the left hand tree is of length $m$ then the right
hand tree is of length $n-m$ where $0<m<n$. The left hand tree will occupy
$2m-1$ levels. There are $\comb{2n-2}{2m-1}$ ways of choosing the levels for
this tree. Hence there are $\comb{2n-2}{2m-1}T(m)$ configurations of the left
hand tree. Moreover, the right hand tree is assigned the remaining levels and
it can be configured $T(n-m)$ ways. Hence there are $\comb{2n-2}{2m-1}T(m)
T(n-m)$ RB trees with $m$ leaves to the left of the root. The result follows
by summing over the different possibilities for $m$.

One should note that the operations corresponding to $\qss $ are not restricted
to consecutive levels. This is because Proposition 2(i) fails if such a
restriction is enforced. Without this property then one would have no means
of constructing a coherent natural isomorphism between any two isomorphic
words composed of the same letters.
The number of RB trees of length $n$ can be calculated by the
recursion formula of Proposition 2(iii). The first six are given in table 2.
\begin{table}[h]
\hspace*{4cm}
\begin{tabular}{|c|c|}
\hline
n & T(n)\\
\hline
1 & 1\\
2 & 2\\
3 & 16\\
4 & 272\\
5 & 7936\\
6 & 353792\\
\hline
\end{tabular}
\caption{The number of RB trees of length $n$.}
\end{table}
The recursive formula for $T(n)$ generates the tangent (or Zag) numbers given
by the $n$th term of the expansion in $\tan x$ thus
\begin{eqnarray}
\tan x & = & \sum_{n=1}^{\infty }T(n)\frac{x^{2n-1}}{(2n-1)!}~.
\end{eqnarray}
Explicitly,
\begin{eqnarray}
T(n) & = & \left. \frac{d^{2n-1}}{dx^{2n-1}}\tan x \right|_{x=0}~.
\end{eqnarray}
Alternatively $T(n)$ is given by the number of up down permutations on $2n-1$
numbers, see Millar {\it et. al.} \cite{jmnsny}.

The underlying $\RBTree $ diagrams for the Dodecagon and Quaddecagon diagrams
are given by figures 13, 14 and 15. The $\qss $--Pentagon is given by figure 5
where the leaf levels occupy the levels above the pivot levels.
\begin{figure}[h]
\hspace*{2cm}
\epsfxsize=300pt
\epsfbox{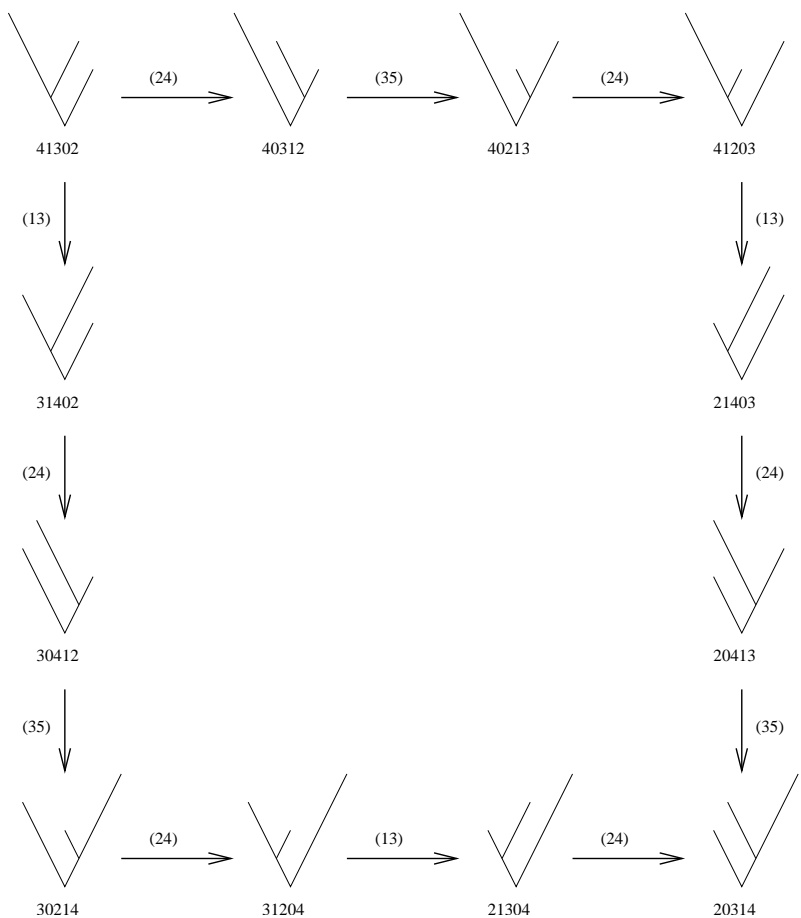}
\caption{Diagram in $\RBTree $ underlying the Dodecagon diagram.}
\end{figure}
\begin{figure}[h]
\hspace*{2cm}
\epsfxsize=300pt
\epsfbox{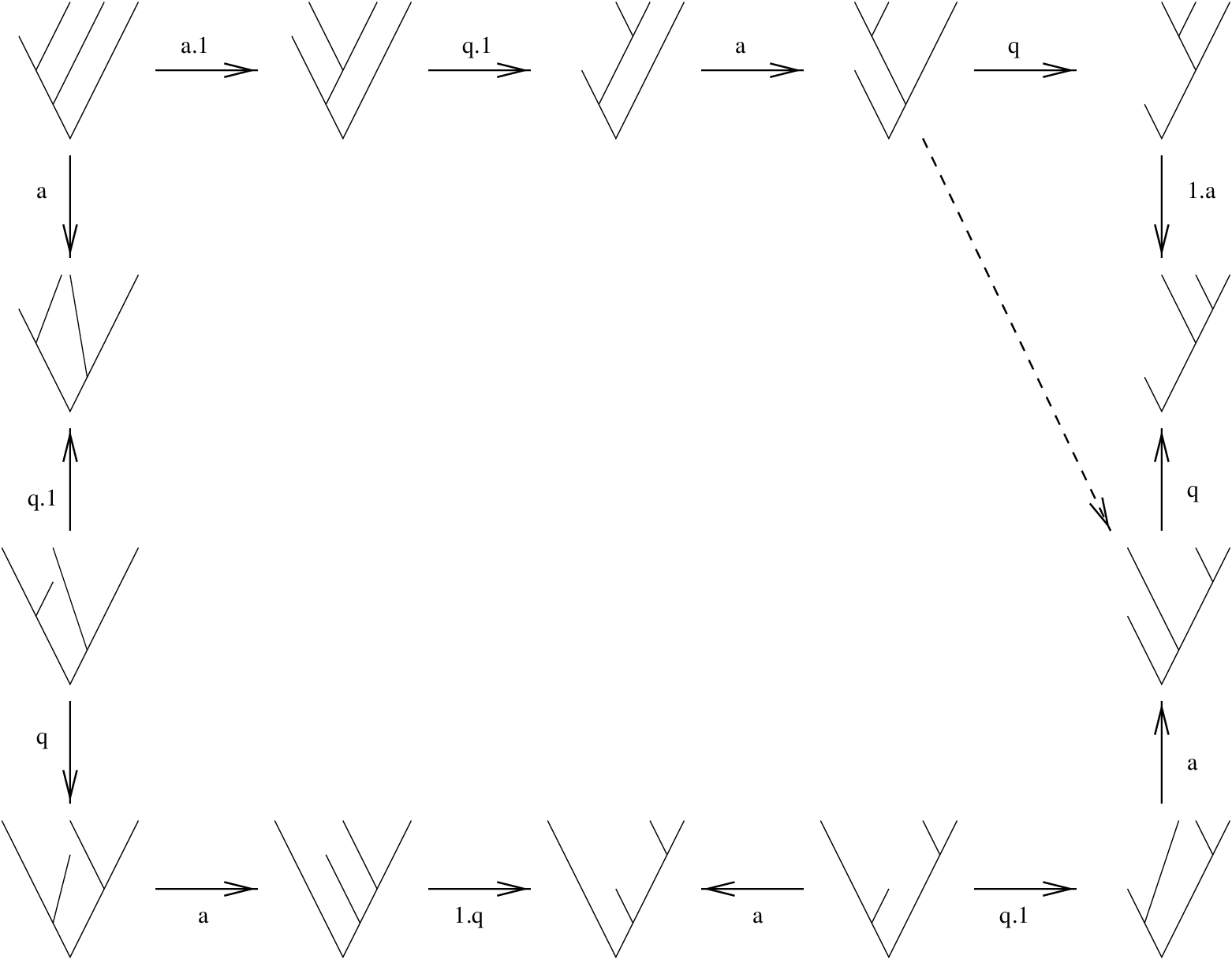}
\caption{Diagram in $\RBTree $ underlying the First Quaddecagon diagram.}
\end{figure}
\begin{figure}[h]
\hspace*{2cm}
\epsfxsize=300pt
\epsfbox{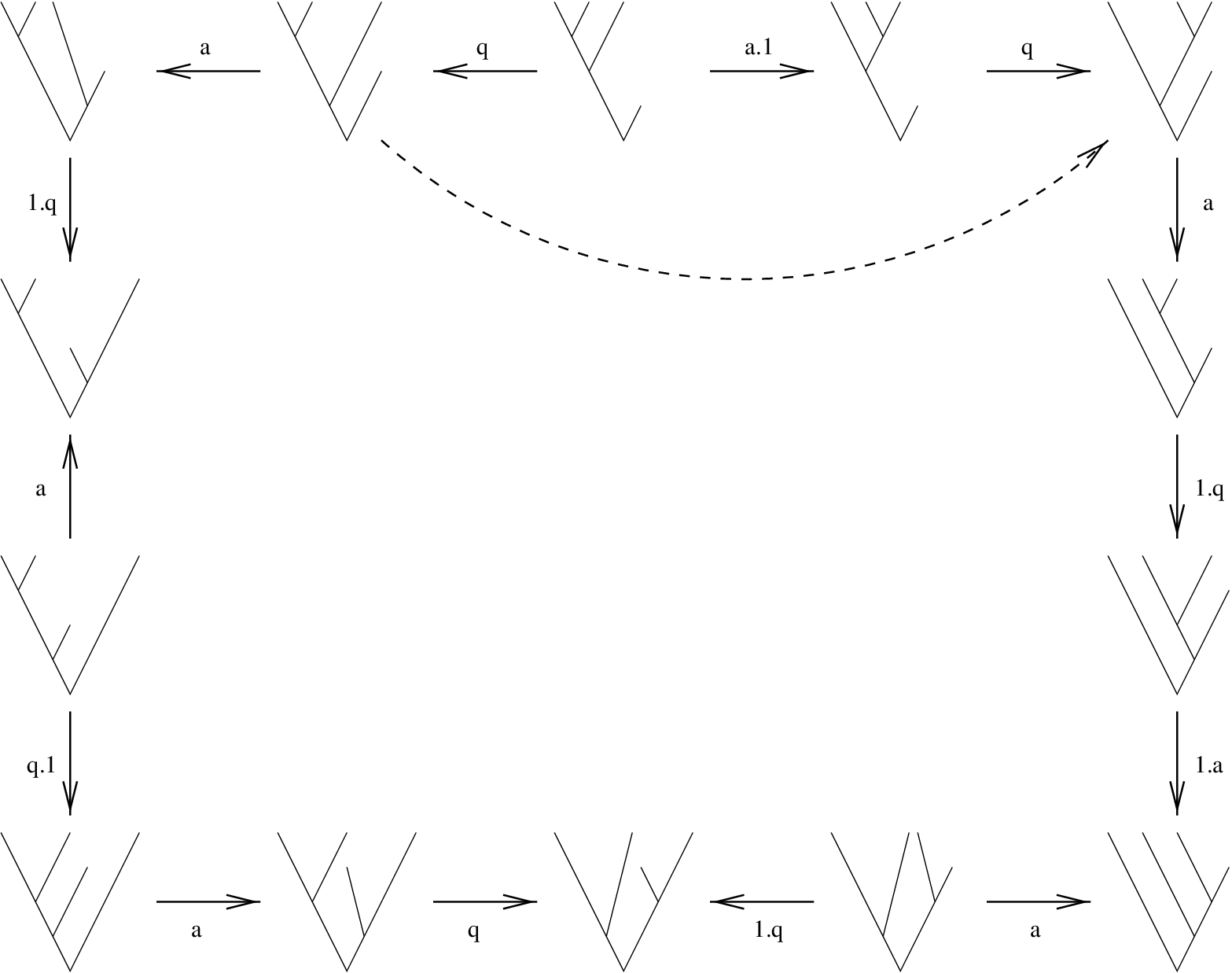}
\caption{Diagram in $\RBTree $ underlying the Second Quaddecagon diagram.}
\end{figure}
Observe that as a direct result of Theorem 1 the Dodecagon and Quaddecagon
diagrams all commute whenever all the leaves in figures 13, 14 and 15 are
replaced by branches.

In general when some of the leaf levels are less than some of the branch
levels constructing a sequence of arrows between any two RB trees is
quite involved. We give an algorithm, called Left Reduction, which constructs
a sequence of primitive arrows reducing any vertex to reduced form. Reduced
form is where all the branch levels are less than all the leaf levels. In fact
it does a little more than this as we shall see. This algorithm plays a central
role in the proof of $\qss $--pseudomonoidal coherence. Given a RB tree,
$B$, the algorithm constructs a sequence of primitive arrows reorganising $B$
into a left--reduced RB tree. We say an RB tree is $l$--left reduced if the
left terminate node of $l$ is not a branch at level $l+1$ and every level
less than $l$ (excluding the root) is a branch and the left terminate for the
next level down. The tree is called left reduced if it is $|B|$--left reduced
or equivalently if the underlying IRB tree is $|B|(|B|-1)...210$. We introduce
some further useful notation. A $k$--cut is where we draw a horizontal line
between the levels $k-1$ and $k$. Every node at a level greater than $k-1$
that is a terminate for a node at a level less than $k$ is called a $k$--cut
node. The line segment between these two nodes cuts the horizontal line. If
the lines in our binary tree do not cross then the $k$--cut nodes have a
unique order along this line. We number them from right to left beginning at
$0$. An example is given in figure 16.
\begin{figure}[h]
\hspace*{2cm}
\epsfxsize=100pt
\epsfbox{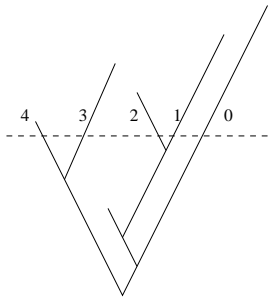}
\caption{A $6$--cut for the RB tree $6480327591(10)$.}
\end{figure}
Note that for an $l$--cut of an $l$--left reduced RB tree the $l$--cut node
positions (except for the highest positioned) are given by the level of the
node for which each is a terminate.

We now give the Left Reduction algorithm for constructing a sequence of
primitive arrows reducing an RB tree, $B$, to a left reduced RB
tree. Note that $\qss $ and $\qss^{-1}$ will not be distinguished. The
levels of the terminates tells you which one is meant. The algorithm
begins with the sequence $B$ which is composed of no arrows. The
algorithm proceeds inductively from the last object of the sequence.\\
\textbf{The Left Reduction Algorithm}
\begin{enumerate}
\item{Locate the largest $l>0$ in the last object of the sequence such that
this object is $l$--left reduced.}
\item{If level $l+1$ is a left terminate of $l$ then (it is a leaf and)
interchange about the pivot $l$.}
\item{Let $r$ be the level with $l+1$ as a terminate. Let $b$ be the branch
level for which the lowest levelled $l$--cut branch node is a terminate. If no
such $b$ exists then go to (x).}
\item{If $r<b$ apply the sequence $\Pi^{b-r-1}_{i=0}\It_{r+i}(\ass^{-1}(1.\qss
)\ass )$.}
\item{If the $l$--cut node attached to $b$ is the left terminate of $l$ then
($b=l$ and) apply $\It_{l}(\qss )$. Thus the $l$--cut node has been lowered to
level $l+1$.}
\item{If $r>b$ apply the sequence $\Pi^{r-b}_{i=1}\It_{r-i}(\ass^{-1}(1.\qss )
\ass )$. Thus the $l$--cut node $b$ has been lowered to level $l+1$.}
\item{If $b=l$ then apply $\It_{l}(\ass^{-1})$ if required to bring into
$(l+1)$--left reduced form and go to (x).}
\item{Apply $\Pi^{l-b}_{j=1}\Pi^{l-b-i}_{i=l-b-j}\It_{b+i}(\ass )$.}
\item{Apply $\Pi^{l-b+1}_{j=1}\Pi^{l-b-i+1}_{i=l-b-j+1}\It_{b+i}(\ass^{-1})$.
Thus the $l$--cut node $b$ has been lowered to level $l+1$.}
\item{If $l<|B|-3$ then go to (i).}
\end{enumerate}
A schematic diagram is given in figure 17 for this algorithm.
\begin{figure}[h]
\hspace*{2cm}
\epsfxsize=250pt
\epsfbox{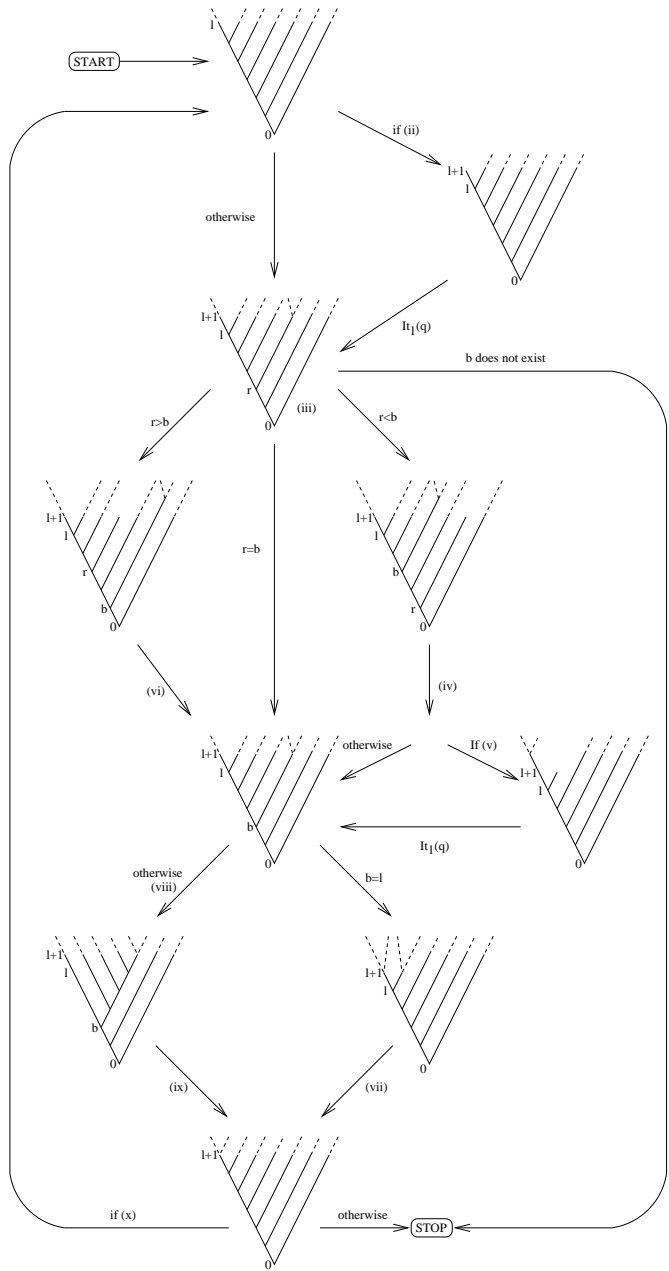}
\caption{A schematic diagram of the Left Reduction algorithm.}
\end{figure}
The following lemma gives us an upper bound on the algorithms complexity.
\begin{lem}
Given an RB tree of length $n$ then the above algorithm reduces it to a left
reduced RB tree using a sequence of $O(n^{2})$ primitive arrows.
\end{lem}
\pf\ Assuming the worst case then for an $l$--cut RB tree we count the
maximum number of primitive arrows: The worst case is the sequence of steps
(i),(iv),(v),(vii),(ix). The maximum number for each is $0$, $3l$, $1$, $l$
and $l+1$ respectively. Summing over $l$ from $1$ to $n-2$ we generate a
sequence of primitive arrows with at most
\begin{eqnarray}
\sum_{l=1}^{n-2}(5l+2) & = & \frac{(5n-1)(n-2)}{2}
\end{eqnarray}
arrows. This is $O(n^{2})$.

We are now in a position to prove a coherence result for $\qss
$--pseudomonoidal categories.
\begin{thm}
  If $\cC $ is a category, $\otimes :\cC \times \cC \rightarrow \cC $ a
  bifunctor and $\ass :\otimes (\otimes \times 1)\rightarrow \otimes
  (1\times \otimes )$ and $\qss :\otimes \rightarrow \otimes $ natural
  isomorphisms then this structure is $\qss $--pseudomonoidal coherent
  if and only if the $\qss $--Pentagon diagram (figure 2), the $\qss
  $--diagrams (figure 6), the Dodecagon diagram (figure 10) and the
  Quaddecagon diagrams (figures 11 and 12) all commute.
\end{thm}
\pf\ The proof is by induction on the length $n$ of RB trees. We define the
rank of a diagram in $\RBTree $ to be the length of any one of its
vertices. The result holds for $n=1,2,3$. Suppose the theorem holds for
diagrams of length $n\geq 3$. Let $D$ be a diagram of length $n+1$. The
induction step is in two parts. Firstly we show that the diagram $D$ is
equivalent to a diagram of the same length whose vertices are all in
reduced form. A vertex is in reduced form if all the branch levels are
below all the leaf levels. A diagram is in reduced form if all of its
vertices are in reduced form. Secondly we show that all diagrams in
reduced form of length $n+1$ commute.

\textbf{Part One:} Let the vertices of $D$ be $a_{0},...,a_{r}$ reading around
the outside and put $a_{r+1}=a_{0}$. We substitute certain arrows that would
cause problems further on in the proof. We replace all arrows $\It_{k}(\ass ),
\It_{k}(\ass^{-1}):a_{i}\rightarrow a_{i+1}$ for which the lowest terminate
level of $k$ is the highest branch of $a_{i}$. We demonstrate how to
substitute for $\It_{k}(\ass )$. The substitution for $\It_{k}(\ass^{-1})$
follows by reversing all arrows. Let $k^{\prime }$ be the left terminate level
of $k$. We assume that $k^{\prime }$ is not the left hand most branch point.
Consider the diagram of figure 22.
\begin{figure}[h]
\epsfxsize=400pt
\epsfbox{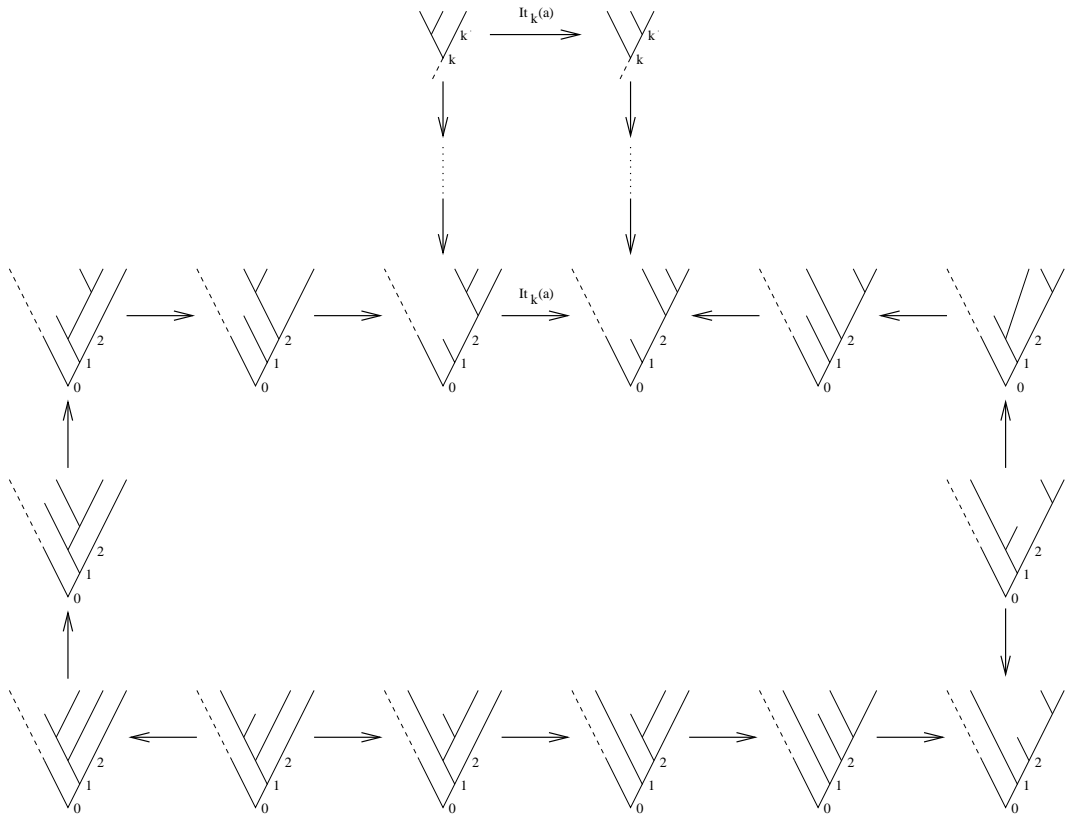}
\caption{Diagram for removing the arrow $\It_{k}(\ass ):a_{i}\rightarrow
a_{i+1}$.}
\end{figure}
We construct two identical sequences of arrows $a_{i}\rightarrow \cdots
\rightarrow b_{i}$ and $a_{i+1}\rightarrow \cdots \rightarrow b_{i+1}$
represented by the vertical sides of the top region. This sequence is
constructed such that the branches $k$ and $k^{\prime }$ and their terminates
remain fixed by each arrow and the highest branch level for each vertex is
$k^{\prime }$. The vertices $b_{i}$ and $b_{i+1}$ are of the form
indicated by diagram. The arrow $b_{i}\rightarrow b_{i+1}$ is $\It_{k}(\ass )$
and the region enclosed commutes by the induction hypothesis and naturality.
We replace the arrow $b_{i}\rightarrow b_{i+1}$ by the fourteen vertex diagram
as indicated. The this region commutes by the first Quaddecagon diagram
(figure 14). Thus we substitute for $\It_{k}(\ass )$ using the perimeter of
the two regions of the diagram in figure 22. If $k^{\prime }$ is the left hand
most branch then a similar substitution is performed using the second
Quaddecagon diagram.

We assume that the diagram $D$ has been substituted as outlined above. For each
$i$ we construct a sequence of arrows $a_{i}\rightarrow \cdots \rightarrow
b_{i}$ using the Left Reduction algorithm. Next we construct a sequence of
arrows $b_{i}\rightarrow \cdots \rightarrow b_{i+1}$ composed of reduced
vertices such that the region enclosed commutes. This will prove the claim of
part one.

Let $k$ be the pivot of $a_{i}\rightarrow a_{i+1}$ then this arrow is either
$\It_{k}(\ass )$ or $\It_{k}(\qss )$. (The inverses are accounted for by
reversing arrows.) Consider the arrow $\It_{k}(\ass ):a_{i}\rightarrow
a_{i+1}$. Let $k^{\prime }$ be the left terminate of $k$. The level
$k^{\prime }$ is not the level of the highest branch. We apply the
Left Reduction algorithm to $a_{i}$ and $a_{i+1}$ until we reach the last
cycle. At this point the last vertices of the sequences constructed coincide.
Moreover, the region enclosed commutes by the induction hypothesis since the
highest branch and its terminates are the same in all vertices. Hence
$b_{i+1}=b_{i}$.

Now consider the arrow $\It_{k}(\qss ):a_{i}\rightarrow a_{i+1}$. Let
$a$ and $b$ be the terminates of $k$ with $a$ to the left of $b$. We allow the
algorithm to cycle through until the last cycle. Let $a_{i}\rightarrow \cdots
\rightarrow c_{i}$ and $a_{i+1}\rightarrow \cdots \rightarrow c_{i+1}$ be the
partially constructed sequences up until the last cycle of the algorithm is
applied. $c_{i}$ and $c_{i+1}$ are $l$--left reduced. The highest branch
either has $a$ and $b$ or neither as terminates.
Let $k^{\prime }$ be the pivot of $a$ and $b$ in $c_{i}$ (and hence $c_{i+1}$).
In the former case the interchange of $a$ and $b$ is natural with respect to
the arrows of the partial sequences. Hence $\It_{k^{\prime }}(\qss ):c_{i}
\rightarrow c_{i+1}$ is a primitive arrow. The region enclosed commutes by
the induction hypothesis. In the latter case neither the levels of the
terminates for highest branch are altered nor is their pivot reattached. Hence
the primitive sequence $\It_{k^{\prime }}(\ass^{-1}(1.\qss )\ass ):c_{i}
\rightarrow c_{i+1}$ encloses a region that is a diagram of length $n$ that
commutes by the induction hypothesis.

It remains to show the commutativity of the region enclosed by the final cycle
of the Left Reduction algorithm. The $c_{i}$ fall into three cases:\\
Case (a): $a\leq l$ or $b\leq l$. One (or both) of $a$, $b$ have been absorbed
and $c_{i+1}=c_{i}$. Hence $b_{i+1}=b_{i}$.\\
Case (b): $a$ and $b$ are terminates of the highest branch. The arrow
$\It_{k^{\prime }}(\qss ):c_{i}\rightarrow c_{i+1}$ commutes with steps (i)
to (viii). Applying these steps we construct sequences $c_{i}\rightarrow \cdots
\rightarrow d_{i}$ and $c_{i+1}\rightarrow \cdots \rightarrow d_{i+1}$. The
arrow $\It_{k^{\prime }}(1.\qss )$ is natural with respect to these sequences.
Moreover, $d_{i}, d_{i+1}$ have highest branch at level $l+1$ with terminates
the leaves $a$ and $b$. The arrow $\It_{l}(\qss ):d_{i}\rightarrow d_{i+1}$
encloses a region that commutes by naturality. Use the final step of the
algorithm to generate the sequences $d_{i}\rightarrow \cdots \rightarrow
b_{i}$ and $d_{i+1}\rightarrow \cdots \rightarrow b_{i+1}$. Hence
these three sequences compose to give a sequence $b_{i}\rightarrow \cdots
\rightarrow d_{i}\rightarrow \cdots \rightarrow d_{i+1}\rightarrow \cdots
\rightarrow b_{i+1}$ where every vertex is in reduced form.\\
Case (c): $a$ and $b$ are attached to different pivots. Note that it is only
possible for one of $a$ and $b$ to be a branch. We continue the algorithm
through to step (vii). Let $c_{i}\rightarrow \cdots \rightarrow d_{i}$ and
$c_{i+1}\rightarrow \cdots \rightarrow d_{i+1}$ be the sequences constructed.
Let $k$ be the pivot of the branch at the $l+1$ level in $d_{i}$ (and hence
$d_{i+1}$). There exists a sequences of primitive arrows $c_{i}\rightarrow
\cdots \rightarrow c_{i+1}$ and $d_{i}\rightarrow \cdots \rightarrow d_{i+1}$
that do not reattach the highest branch nor interchange its terminates. Hence
the regions $a_{i}a_{i+1}c_{i+1}c_{i}$ and $c_{i}c_{i+1}d_{i+1}d_{i}$ commute
by the induction hypothesis and naturality. Use the final steps of the
algorithm to generate the sequences $d_{i}\rightarrow \cdots \rightarrow
b_{i}$ and $d_{i+1}\rightarrow \cdots \rightarrow b_{i+1}$. Hence these three
sequences compose to give a sequence $b_{i}\rightarrow \cdots \rightarrow
d_{i}\rightarrow \cdots \rightarrow d_{i+1}\rightarrow \cdots \rightarrow
b_{i+1}$ where every vertex is in reduced form.

\textbf{Part Two:}
We now prove the second step which is that a diagram $D$ of length $n+1\geq 5$
whose vertices are in reduced form commutes. We label the leaf levels using
an $(n+1)$--cut. The leaf $v$ is assigned the number $2n-1-l(v)$. This gives
$\qss $ the action of a permutation on the leaf levels. Since $D$ is in
reduced form we are restricted to those leaves that interchange branch levels
and those that interchange leaf levels. the former we call terminal arrows and
the latter internal arrows. The latter we will show allow us to construct
generators for the permutation group $S_{n+1}$. If $D$ contains only internal
arrows then $D$ commutes by Theorem 1.

Suppose that $D$ contains at least one terminal arrow. We partition $D$ into
sequences of arrows such that each sequence has exactly one arrow for
interchange of leaf levels. Let a typical partition sequence be
$a_{m}\rightarrow \cdots \rightarrow a_{p}$. We will show that any alternative
sequence $b_{0},...,b_{q}$ with $b_{0}=a_{m}$ and $b_{q}=a_{p}$ containing
only one interchange of leaf levels encloses a region that commutes. Suppose
$a_{k}\rightarrow a_{k+1}$ and $b_{l}\rightarrow b_{l+1}$ are the interchanges
of leaf levels. The levels they interchange are identical and occur in the
same adjacent positions. Hence there are sequences of internal arrows from
$a_{k}$ to $b_{l}$ and $a_{k+1}$ to $b_{l+1}$ such that the region enclosed
commutes by naturality. This is depicted in figure 19.
\begin{figure}[h]
\hspace*{2cm}
\epsfxsize=200pt
\epsfbox{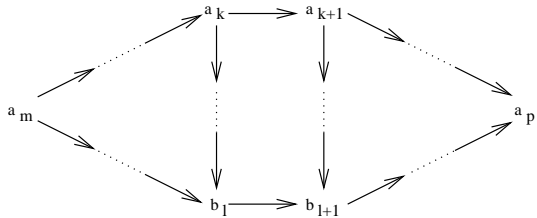}
\caption{Diagram showing $\tau_{i}:a_{m}\rightarrow a_{p}$ is well--defined.}
\end{figure}
Moreover, the regions $a_{m}a_{k}b_{l}b_{0}$ and $a_{p}a_{k+1}b_{l+1}b_{q}$
commute. Hence the two sequences between $a_{m}$ and $a_{p}$ compose to give
the same natural isomorphism.

We define the adjacent transpositions between two objects $a$ and $b$ to be
the arrows $\tau_{i}=(i\ i+1):a\rightarrow b$ interchanging the leaf levels in
positions $i$ and $i+1$ that is given by any sequence of arrows between
the objects $a$ and $b$ with only one terminal arrow. This definition is
well--defined by the previous paragraph. Note that each adjacent transposition
is a family of arrows index by source and target. We partition $D$ into
sections where each section is a sequence of arrows containing only one
terminal arrow. Let $a_{0},...,a_{r}$ be the boundary vertices of these
sections with $a_{r+1}=a_{0}$. Each sequence $a_{k}\rightarrow \cdots
\rightarrow a_{k+1}$ of $D$ corresponds to an adjacent transposition
$\tau_{i}:a_{k}\rightarrow a_{k+1}$. We now show that the adjacent
transpositions satisfy
\begin{eqnarray}
\tau_{i}^{2}=1~, & & i=1,...,n \\
\tau_{i}\tau_{j}=\tau_{j}\tau_{i}~, & & 1\leq i<j-1\leq n-1\\
\tau_{i}\tau_{i+1}\tau_{i}=\tau_{i+1}\tau_{i}\tau_{i+1}~, & \ph{space} &
i=1,...,n-2
\end{eqnarray}
These are the generating relations for $S_{n}$ from which it follows that
$D$ commutes.\\
Property (i): This follows because any sequence defining $\tau_{i}:a
\rightarrow b$ gives $\tau_{i}:b\rightarrow a$ by inverting the sequence.\\
Property (ii): Given $\tau_{i}:a\rightarrow c$, $\tau_{j}:c\rightarrow b$,
$\tau_{j}:a\rightarrow d$ and $\tau_{i}:d\rightarrow b$. We construct
sequences of internal arrows from $a,b,c,d$ to vertices in $(n-3)$--left
reduced form as given by the diagram in figure 20.
\begin{figure}[h]
\hspace*{2cm}
\epsfxsize=300pt
\epsfbox{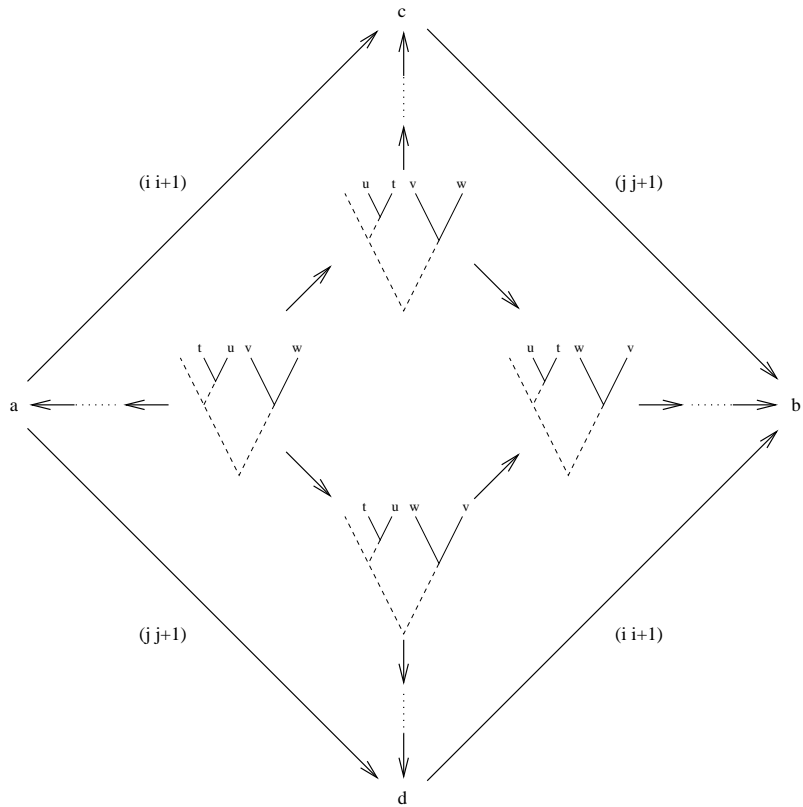}
\caption{The property $\tau_{i}\tau_{j}=\tau_{j}\tau_{i}$ whenever $1\leq
i<j-1\leq n-1$.}
\end{figure}
The branches shown are at the $n-1$ and $n-2$ levels and the terminates are in
the $i,i+1,j,j+1$ positions with their levels labelled $t,u,v,w$. The centre
region is a natural square. All other regions commute because the definition
of $\tau_{i}, \tau_{j}$ is well--defined.\\
Property (iii): Given $\tau_{i}:a\rightarrow c$, $\tau_{i+1}:c\rightarrow d$,
$\tau_{i}:d\rightarrow b$, $\tau_{i+1}:a\rightarrow e$, $\tau_{i}:e\rightarrow
f$ and $\tau_{i+1}:f\rightarrow b$; we construct sequences of arrows as given
by the diagram in figure 21.
\begin{figure}[h]
\epsfxsize=400pt
\epsfbox{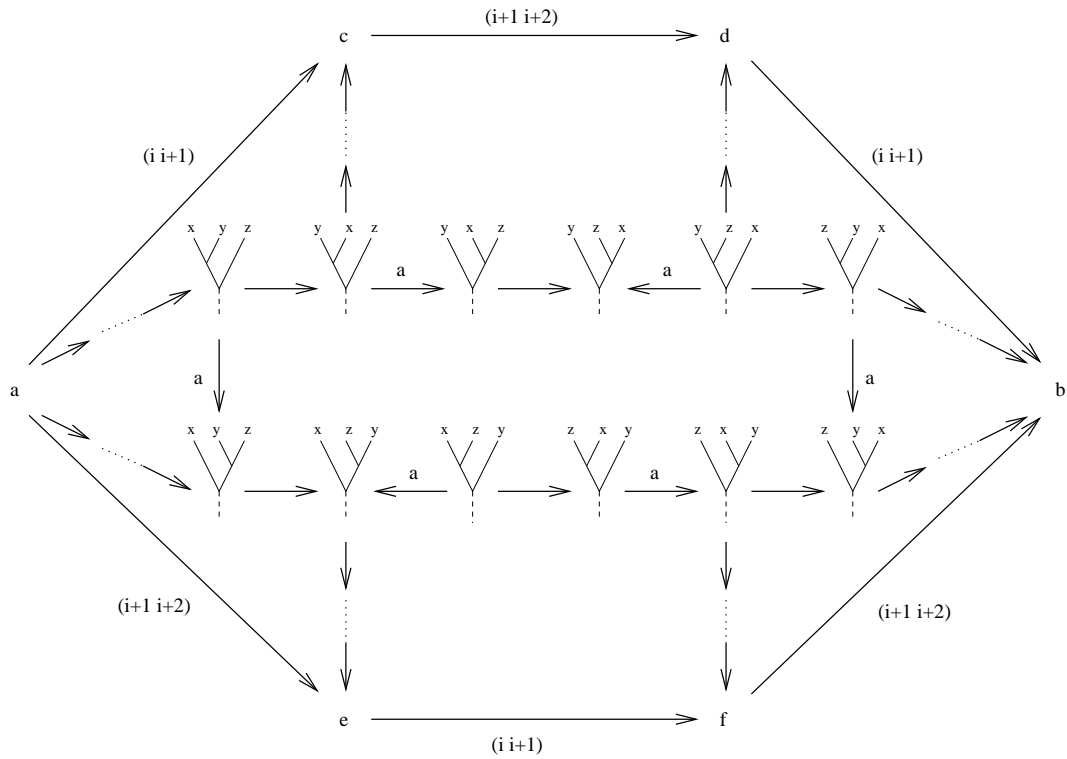}
\caption{The property $\tau_{i}\tau_{i+1}\tau_{i}=\tau_{i+1}\tau_{i}
\tau_{i+1}$ whenever $1\leq i\leq n$.}
\end{figure}
The central diagram is the Dodecagon diagram (figure 13) and commutes.
The vertices of this diagram are taken to be in $(n-3)$--left reduced
form. The branches shown are at the $n-1$ and $n-2$ levels and the
terminates are in the $i,i+1,i+2$ positions with their levels labelled
$x,y,z$. The sequences from $a,b,c,d,e,f$ to the vertices of the centre
region are taken to be composed only of internal arrows. Moreover, all
the regions enclosed around the centre region commute because
$\tau_{i},\tau_{i+1}$ are well--defined. This completes the proof of
Theorem 2.

\section{$\qss $--Braided Pseudomonoidal Categories}             %

The step from pseudomonoidal categories to $\qss $--pseudomonoidal
categories was a natural one. Unfortunately one is faced with two
fourteen vertex diagrams. Closer inspection of the diagram conditions
shows that the $\qss $--Pentagon diagram concerns the action of $\qss $
with branch terminates; the Dodecagon diagram concerns the action of
$\qss $ with leaf terminates; and the Quaddecagon and Dodecagon diagrams
concern the action of $\qss $ with mixed--typed terminates. If we
consider the full subcategory, $\RRBTree $ of reduced resolved binary
trees, of $\RBTree $ then the two Quaddecagon diagrams are redundant.
This leads to a coherence result that is weaker than $\qss
$--pseudomonoidal coherence admitting a braid structure. The Dodecagon
diagram is the Yang--Baxter equation for this braid. We begin by
weakening our notion of $\qss $--premonoidal category.
\begin{defn}
A $\qss $--braided pseudomonoidal category is a quadtruple $(\cC ,\otimes ,\ass
,\qss )$ where $(\cC ,\otimes ,\ass )$ is a pseudomonoidal category and $\qss :
\otimes \rightarrow \otimes $, called the $\qss $--braid, is a natural
automorphism satisfying the Dodecagon diagram (figure 10).
\end{defn}
The $\qss $ arising from the pseudomonoidal structure is in general
different from the $\qss $--braid. The first is a natural automorphism
for the functor $\otimes (\otimes \times \otimes ):\cC^{4}\rightarrow \cC $
corresponding to branch level interchange. The second is a natural
automorphism for the functor $\otimes :\cC^{2}\rightarrow \cC $ corresponding
to leaf level interchange. Every $\qss $--pseudomonoidal category is a
$\qss $-braided strong pseudomonoidal category. For $\qss $--pseudomonoidal
categories the $\qss $ and $\qss $--braid are identified as the same
natural automorphism.

We collect some results about $\RRBTree $. A reduced RBtree of length $n$ is
uniquely represented by an ordered pair of linear orderings. The first linear
ordering of the numbers $0,1,...,n-2$ represents the underlying IRB tree. The
second linear ordering of the numbers $1,2,...,n$ is assigned to each leaf
reading from left to right across the tree. The leaf level is given by adding
$n-2$ to the number assigned by the second linear ordering. An example is
given in figure 22.
\begin{figure}[h]
\hspace*{4cm}
\epsfxsize=60pt
\epsfbox{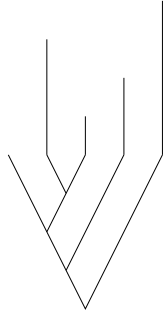}
\caption{The reduced resolved binary tree $(2310,14235)$}
\end{figure}
We call the first sequence ($2310$) the branch structure and the second
sequence ($14235$) the leaf structure.
\begin{prop}$\ph{nothing}$
\begin{enumerate}
\item{For any two RRB trees of length $n$ there is a finite sequence of
primitive arrows transforming one into the other.}
\item{Every RRB tree of length $n$ is the source of at most $n-1$ distinct
primitive arrows.}
\item{There are $n!(n-1)!$ RRB trees of length $n$.}
\end{enumerate}
\end{prop}
\pf\ (i) Let $B$ and $B^{\prime }$ be two RRB trees of length $n$. Every
permutation is the product of a finite number of transpositions. Hence the
result follows if very transposition of $B$ to $B^{\prime }$ can be
constructed from a sequence of primitive arrows. We construct (by Proposition
1) a sequence of internal arrows from $B$ to an RRB tree with the $i$ and
$i+1$ terminates of the branch at level $n-1$. Next we interchange $i$ and
$i+1$. Finally we construct a sequence of internal arrows to $B^{\prime }$.\\
(ii) We prove by induction on the length of the RRB tree. Suppose the result
holds for RRB trees of length $n$. Given an RRB tree of length $n+1$ we
consider the highest branch level to be a leaf. Then there are at most $n-1$
distinct primitive arrows. Including the highest branch we obtain an extra
arrow for interchange and possibility an arrow for reattachment. The second
arrow only occurs if the branch is joined to a pivot whose other terminate
is a leaf. In this situation we over--counted an interchange of
leaf levels. Hence there are at most $n$ primitive operations.\\
(iii) There are $(n-1)!$ of arranging the branch structure and $n!$ ways of
rearranging the leaf structure.

The notion of RRB coherence follows directly from Theorem 2, part two of the
proof. However, we can obtain a more general braided coherence result. We
construct from $\RRBTree $ a braid category denoted $\qB $. The objects are
RRB trees. The arrows are reattachment and interchange of branches as for
$\RRBTree $. The $\qss :\otimes \rightarrow \otimes $ natural automorphism
corresponds to primitive arrows represented pictorially by the identification
with a braid as given in figure 23.
\begin{figure}[h]
\hspace*{2cm}
\epsfxsize=150pt
\epsfbox{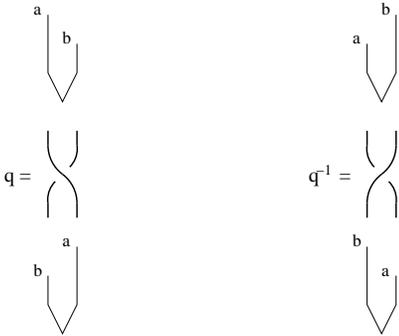}
\caption{The $\qss $--braid primitive arrows.}
\end{figure}
These arrows (or $\qss $--braids) act without regard for the level of the
leaves. That is, there is no restriction that $a>b$ for $\qss $ and $b<a$
for $\qss^{-1}$. The primitive $\qss $--braids are arrows that interchange
two adjacent leaf levels provided the leaf levels interchanged are terminates
of the same branch. The Dodecagon diagram arises from the $\qss $--Yang--Baxter
diagram given in figure 24.
\begin{figure}[h]
\hspace*{2cm}
\epsfxsize=150pt
\epsfbox{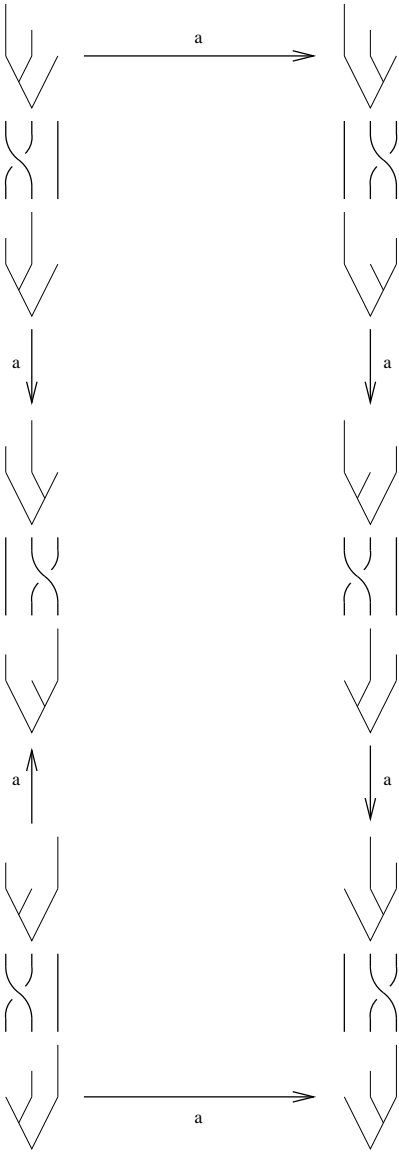}
\caption{The $\qss $--Yang--Baxter diagram.}
\end{figure}
There is a lot of freedom at the branch level. One way to avoid this is to
choose an RRB tree $a_{n}$ for each $\qB_{n}$ where $\qB_{n}$ is the full
subcategory of $\qB $ containing all the RRB trees of length $n$. We call any
set of objects $( a_{n} )_{n\in \bN }$ where $a_{n}$ is of length $n$ a frame
for $\qB $. Let $\pi \in S_{n}$ and $b_{1}...b_{n}$ be the leaf structure for
$a_{n}$. We define $\pi a_{n}$ to be the RRB tree with the same branch
structure as $a_{n}$ and leaf structure $b_{\pi (1)}...b_{\pi (n)}$ given by
permuting the leaf structure of $a_{n}$. The primitive $\qss $--braids of
length $n$, denoted $\sigma_{i}$ where $i=1,...,n$; in the frame
$(a_{n})_{n\in \bN }$ are sequences of primitive arrows $\sigma_{i}:a_{n}
\rightarrow \cdots \rightarrow (i\ i+1)a_{n}$, where $1\leq i\leq n$.
Moreover, these sequences can be any choice constructed using the rules of
$\RRBTree $. The Dodecagon diagram implies the braid generating condition
given in figure 25.
\begin{figure}[h]
\hspace*{2cm}
\epsfxsize=80pt
\epsfbox{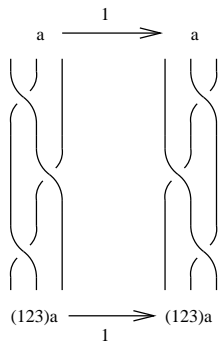}
\caption{The braid generating condition.}
\end{figure}
The primitive $\qss $--braids generate the Artin braid group. Under the
functor $\can $ this $\qss $--braid group gives a braid coherence result.
More generally we have the following coherence result.
\begin{thm}
  If $\cC $ is a category, $\otimes :\cC \times \cC \rightarrow \cC $ a
  bifunctor and $\ass :\otimes (\otimes \times 1)\rightarrow \otimes
  (1\times \otimes )$ and $\qss :\otimes \rightarrow \otimes $ natural
  isomorphisms then this structure is $\qss $--braided pseudomonoidal
  coherent if and only if the $\qss $--Pentagon diagram (figure 2), the
  $\qss $--Square diagrams (figure 6) and the Dodecagon diagram (figure
  10) all commute.
\end{thm}
\pf\ By part two of the proof of Theorem 2.

\section{Unital Pseudomonoidal Categories}                       %

We now consider the issue of including an identity into pseudomonoidal
and $\qss $--pseudomonoidal structures. This requires an identity object
and natural isomorphisms for contracting the identity object on the left
and right. The result is a unital pseudomonoidal or unital $\qss
$--pseudomonoidal category. We only consider the pseudomonoidal case in what
follows because the $\qss $--pseudomonoidal case follows similarly. The
unital pseudomonoidal structure is monoidal when $\qss =1$.  Moreover,
the coherent diagrams for a restricted premonoidal structure are
monoidal whenever the identity object is involved. This is some what
dissatisfying and will be resolved in the next section.
\begin{defn}
A unital pseudomonoidal category is a sextet $(\cC ,\otimes ,\ass ,\lid ,\rid
,e)$ such that $(\cC ,\otimes ,\ass )$ is a pseudomonoidal category (and hence
allows construction of the natural automorphism $\qss $), $e$ is an
object of $\cC $ called the identity object and $\lid :e\otimes \ubar
\rightarrow 1$ and $\rid :\ubar \otimes e\rightarrow 1$ are natural
isomorphisms satisfying the Triangle diagram (figure 26).
\end{defn}
\begin{figure}[h]
$\begin{diagram}
\putVtriangle<1`1`1;600>(600,0)[(\alpha \otimes e)\otimes \beta `\alpha
\otimes (e\otimes \beta )`\alpha \otimes \beta ;\ass_{\alpha ,e,\beta }`
\rid_{\alpha }\otimes 1_{\beta }`1_{\alpha }\otimes \lid_{\beta }]
\end{diagram}$
\caption{The Triangle diagram ($T$).}
\end{figure}
\begin{figure}[h]
$\begin{diagram}
\putVtriangle<1`1`1;600>(600,0)[(e\otimes \alpha )\otimes \beta `e\otimes
(\alpha \otimes \beta )`
\alpha \otimes \beta ;\ass_{e,\alpha ,\beta }`\lid_{\alpha }\otimes 1_{\beta
}`\lid_{\alpha \otimes \beta }]
\putVtriangle<1`1`1;600>(2400,0)[(\alpha \otimes \beta )\otimes e`\alpha
\otimes (\beta \otimes e)`\alpha \otimes \beta ;\ass_{\alpha ,\beta ,e}`
\rid_{\alpha \otimes \beta }`1_{\alpha }\otimes \rid_{\beta }]
\end{diagram}$
\caption{The redundant Triangle diagrams ($T_{l}$ and $T_{r}$ respectively).}
\end{figure}
We assume for now that all the triangle diagrams hold (figures 26 and 27).
We will shortly show that $\qss $ is the identity whenever $e$ is an index.
Whence from the monoidal case only the Triangle diagram of figure 26 is
required. The redundancy amongst the triangle diagrams was first pointed out
by Kelly \cite{gk}.
\begin{figure}[h]
$\begin{diagram}
\putVtriangle<1`1`1;600>(300,0)[(e\otimes \alpha )\otimes (\beta \otimes
\gamma)`(e\otimes \alpha )\otimes (\beta \otimes \gamma )`\alpha \otimes (\beta
\otimes \gamma );\qss_{e,\alpha ,\beta ,\gamma }`\lid_{\alpha }\otimes 1_{
\beta \otimes \gamma }`\lid_{\alpha }\otimes 1_{\beta \otimes \gamma }]
\putVtriangle<1`1`1;600>(2400,0)[(\alpha \otimes e)\otimes (\beta \otimes
\gamma )`(\alpha \otimes e)\otimes (\beta \otimes \gamma )`\alpha \otimes
(\beta \otimes \gamma );\qss_{\alpha ,e,\beta ,\gamma }`\rid_{\alpha }\otimes
1_{\beta \otimes \gamma }`\rid_{\alpha }\otimes 1_{\alpha \otimes \gamma }]
\putVtriangle<1`1`1;600>(300,-900)[(\alpha \otimes \beta )\otimes (e\otimes
\gamma )`(\alpha \otimes \beta )\otimes (e\otimes \gamma )`(\alpha \otimes
\beta )\otimes \gamma ;\qss_{\alpha ,\beta ,e,\gamma }`1_{\alpha \otimes \beta
}\otimes \lid_{\gamma }`1_{\alpha \otimes \beta }\otimes \lid_{\gamma }]
\putVtriangle<1`1`1;600>(2400,-900)[(\alpha \otimes \beta )\otimes (\gamma
\otimes e)`(\alpha \otimes \beta )\otimes (\gamma \otimes e)`(\alpha \otimes
\beta )\otimes \gamma ;\qss_{\alpha ,\beta ,\gamma ,e}`1_{\alpha \otimes \beta
}\otimes \rid_{\gamma }`1_{\alpha \otimes \beta }\otimes \rid_{\gamma }]
\end{diagram}$
\caption{The triangle diagrams $Q_{1}$, $Q_{2}$, $Q_{3}$ and $Q_{4}$ satisfied
by $\qss $.}
\end{figure}
The $\qss $ natural automorphism satisfies the triangular diagrams of figure 28
because of the following lemma.
\begin{lem}
Let $\cC $ be a category, $\otimes :\cC \times \cC \rightarrow \cC $ a
bifunctor, $e$ an object of $\cC $, $\ass $, $\qss $, $\lid $ and $\rid $
natural isomorphisms corresponding to associativity, deformation, left identity
and right identity respectively. The following equivalences hold.
\begin{enumerate}
\item{If $T_{l}$ commutes then the $\qss $--Pentagon diagram commutes if and
only if $Q_{1}$ commutes.}
\item{If $T_{r}$ commutes then the $\qss $--Pentagon diagram commutes if and
only if $Q_{4}$ commutes.}
\item{If $T$ commutes then the commutativity of any two of the
$\qss $--Pentagon diagram, $T_{l}$ and $Q_{2}$ implies the third.}
\item{If $T$ commutes then the commutativity of any two of the
$\qss $--Pentagon diagram, $T_{r}$ and $Q_{3}$ implies the third.}
\end{enumerate}
\end{lem}
\pf\ Each part follows by considering the $\qss $--Pentagon diagram with
$a=e$, $d=e$, $b=e$ and $c=e$ respectively.

The bottom two arrows for each diagram in figure 26 compose to give an identity
arrow. Hence we have the following proposition.
\begin{prop}
Let $(\cC ,\otimes ,\ass ,\lid ,\rid ,e)$ be a unital pseudomonoidal category
then $\qss_{\alpha ,\beta ,\gamma ,\delta }=1_{(\alpha \otimes \beta )\otimes
(\gamma \otimes \delta )}$ whenever any one of $\alpha ,\beta ,\gamma ,\delta $
is $e$.
\end{prop}

We move onto the issue of coherence. The coherence structure of interest is
the category of internally resolved binary trees with nodules. A nodule is an
open circle which may be attached to any leaf. Note that each nodule is
uniquely determined by the adjacent internal node levels in the IRB tree number
sequence. Whence a nodule is represented by placing a dot between these two
levels. An example is given in figure 29.
\begin{figure}[h]
\hspace*{2cm}
\epsfxsize=100pt
\epsfbox{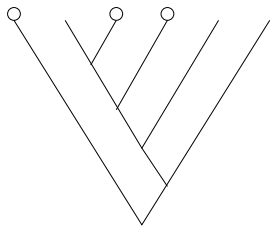}
\caption{The IRB tree with nodules described by $\cdot 04\cdot 3\cdot 21$.}
\end{figure}
The length $|B|$ of an IRB tree $B$ with nodules is defined to be the number
of leaves without nodules. The number of nodules is given by
$\nodes{B}=\hght{B}-|B|+1$. The groupoid of IRB trees with nodules is
denoted $\IRNBTree $. The objects are IRB trees with nodules. The
primitive arrows are inherited from $\IRBTree $ with the addition of
primitive arrows for pruning and grafting nodules. The primitive arrows
for the pruning of a nodule are given in figure 30.
\begin{figure}[h]
\hspace*{2cm}
\epsfxsize=100pt
\epsfbox{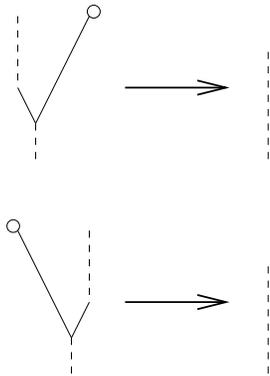}
\caption{The primitive arrows for pruning nodules corresponding to iterates of
left and right identity natural isomorphisms.}
\end{figure}
The dashed lines represent attachment sites to the remaining edges and
nodes of the binary tree. Note that the source of the pruning arrow
requires the level of a branch terminate to be on the next level up. The
pruning arrow deletes a nodule from the tree and reduces by one the level of
all nodes greater than the pivot level. The primitive arrows for grafting are
given by the inverse of pruning. We have the following nodule version of
Proposition 1.
\begin{prop}$\ph{nothing}$
\begin{enumerate}
\item{Given two IRB trees with nodules of length $n$ then there is a
finite sequence of primitive arrows transforming one into the other.}
\item{Every IRB tree of height $n$ with $m$ nodules is the source of no
more than $5n-1+m$ distinct primitive arrows.}
\item{There are $\comb{n+1}{m}n!$ IRB trees of height $n$ with $m$ nodules.}
\end{enumerate}
\end{prop}
\pf\ (i) This follows from Proposition 1 since every IRB tree with
nodules can have all of its nodules pruned and nodules may be grafted onto any
edge.\\
(ii) We have that there are no more than $n-1$ $\IRBTree $ primitive arrows.
There are $m$ primitive arrows removing nodules and $2|E|$ adding nodules
since a nodule may be attached to either side of every edge. Finally the
number of edges is $|E|=2n$ so summing the possibilities we obtain the
result.\\
(iii) There are $n!$ IRB trees of height $n$ and for each such tree
there are $\comb{n+1}{m}$ ways of filling the $n+1$ leaves with
$m$ nodules.

The objects $B$ of $\IRNBTree $ give rise to functors as in the nodule free
situation with the following differences. The nodules stand in for the object
$e$. If we let $\Endo_{\cC }(e)$ denote the full subcategory of $\cC $
with the single object $e$ then the functor $\can (B):\Endo_{\cC }(e)^{
\nodes{B}}\times \cC^{|B|}\rightarrow \cC $ is given by contracting objects and
arrows according to $B$. For example the tree of figure 27 gives a functor
taking $(c_{1},c_{2},c_{3})\mapsto e\otimes ((((c_{1}\otimes e)\otimes
e)\otimes c_{2})\otimes c_{3})$ and $(f_{1},f_{2},f_{3},g_{1},g_{2}, g_{3})
\mapsto f_{1}\otimes (((g_{1}\otimes f_{2})\otimes f_{3})\otimes g_{2})
\otimes g_{3})$ where $c_{i}$ are objects in $\cC $, $f_{i}$ are
arrows in $\Endo_{\cC }(e)$ and $g_{i}$ are arrows in $\cC $. We have the
following expected coherence result.
\begin{thm}
A pseudomonoidal category $(\cC ,\otimes ,\ass )$ with identity object $e$ and
natural isomorphisms $\lid :e\otimes \ubar \rightarrow 1$ and
$\rid :\ubar \otimes e\rightarrow 1$ is unital pseudomonoidal coherent if and
only if the Triangle diagram (figure 26) commutes.
\end{thm}
\pf\ Let $D$ be a diagram in $\cD (\IRNBTree )$. We define the rank of each
vertex to be the number of nodules that it contains. The rank of $D$ is
defined to be the maximum of vertex ranks. Proof is by induction on diagram
rank. The result holds for all diagrams of rank $0$ by Theorem 1. Suppose the
result is true for diagrams of rank $n+1$. Let
$a_{0},...,a_{r}$ be the vertices for some diagram $D$ of rank $n+2$ given by
reading around the outside. We identify $a_{r+1}$ with $a_{0}$. We connect
each vertex $a_{k}$ of $D$ to a vertex $b_{k}$ obtained as follows: If $a_{k}$
has rank $n+1$ then remove the left hand most nodule using the primitive
arrow for its removal. If the rank of $a_{k}$ is less than $n+2$
then define $b_{k}=a_{k}$. Connect $b_{k}$ to $b_{k+1}$ using a
primitive arrow to form the square $a_{k}a_{k+1}b_{k+1}b_{k}$. This is checked
by considering all possible cases. The squares obtained are natural or
correspond to the triangle diagrams of figures 26, 27 and 28. Finally the
diagram with vertices $b_{0},...,b_{r}$ constructed is of rank less than or
equal to $n$ and by the induction hypothesis commutes. Hence the diagram $D$
commutes.

\section{Symmetric $\qss $--Pseudomonoidal Categories}           %

We now weave a natural isomorphism for commutativity into the picture.
We consider a commutativity natural isomorphism that is symmetric.
Moreover, additional symmetric natural isomorphisms can be constructed
from $\com $ and $\qss $. The weakening of a symmetric $\qss
$--pseudomonoidal structure to a symmetric $\qss $--braided pseudomonoidal
structure is an obvious step. The $\qss $--braid is symmetric. It is
only this latter structure that admits what should properly be called a
symmetric $\qss $--monoidal category.

Define the flip functor to be $\tau :\cC \times \cC \rightarrow \cC \times
\cC $ where $(a,b)\mapsto (b,a)$.
\begin{defn}
A symmetric $\qss $--pseudomonoidal category is a pentuple $(\cC ,\otimes
,\ass ,\qss ,\com )$ where $(\cC ,\otimes ,\ass  ,\qss )$ is a
$\qss $--pseudomonoidal category, $\com :\otimes \rightarrow \otimes \tau $ is
a natural isomorphism such that $\com $ and $\qss $ are symmetric and the
Square diagram (figure 31), the Hexagon diagram (figure 32) and the Decagon
diagrams (figures 33 and 34) all commute.
\end{defn}
\begin{figure}[h]
$\begin{diagram}
\putsquare<1`1`-1`-1;600`600>(600,0)[\alpha \otimes \beta `\beta \otimes
\alpha `\alpha \otimes \beta `\beta \otimes \alpha ;\com_{\alpha ,\beta }`
\qss_{\alpha ,\beta }`\qss_{\beta ,\alpha }`\com_{\beta ,\alpha }]
\end{diagram}$
\caption{The Square diagram.}
\end{figure}
\begin{figure}[h]
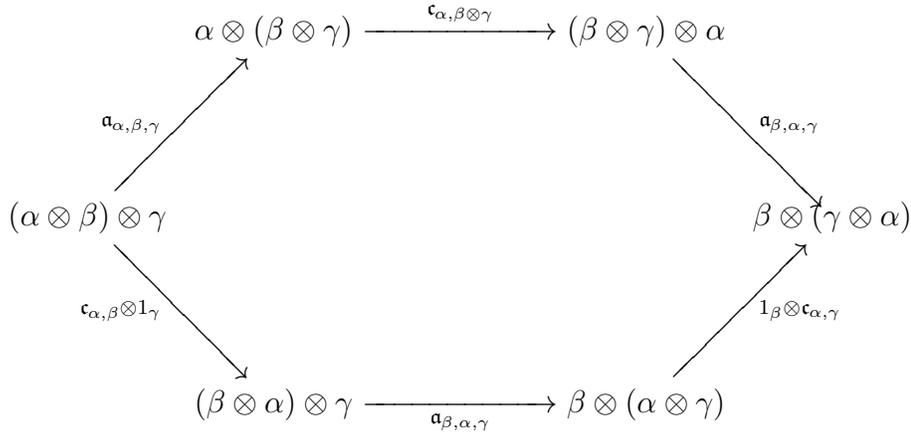

$\begin{diagram}
\putdtriangle<-1`0`0;600>(600,0)[\alpha \otimes (\beta \otimes \gamma )`
(\alpha \otimes \beta )\otimes \gamma `;\ass_{\alpha ,\beta ,\gamma }``]
\putqtriangle<0`1`0;600>(600,-600)[\ph{(\alpha \otimes \beta )\otimes \gamma
}``(\beta \otimes \alpha )\otimes \gamma ;`\com_{\alpha ,\beta }\otimes
1_{\gamma }`]
\puthmorphism(1200,600)[\ph{\alpha \otimes (\beta \otimes \gamma )}`(\beta
\otimes \gamma )\otimes \alpha `\com_{\alpha ,\beta \otimes \gamma }]{1200}{1}a
\puthmorphism(1200,-600)[\ph{(\beta \otimes \alpha )\otimes \gamma }`\beta
\otimes (\alpha \otimes \gamma )`\ass_{\beta ,\alpha ,\gamma }]{1200}{1}b
\putbtriangle<0`1`0;600>(2400,0)[\ph{(\beta \otimes \gamma )\otimes \alpha )}
``\beta \otimes (\gamma \otimes \alpha );`\ass_{\beta ,\alpha ,\gamma }`]
\putptriangle<0`0`-1;600>(2400,-600)[`\ph{\beta \otimes (\gamma \otimes \alpha
)}`\ph{\beta \otimes (\alpha \otimes \gamma )};``1_{\beta }\otimes
\com_{\alpha ,\gamma }]
\end{diagram}$
\caption{The Hexagon diagram.}
\end{figure}
\begin{figure}[h]
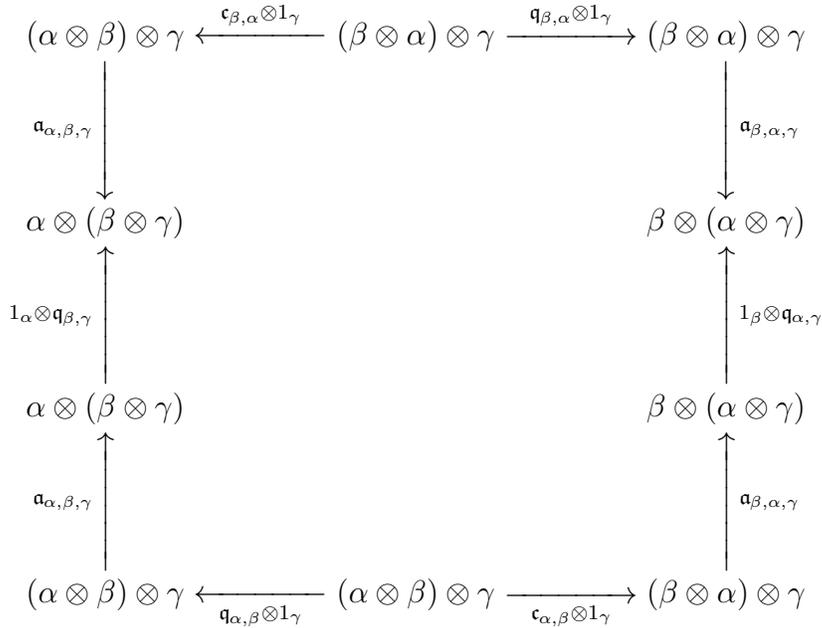

$\begin{diagram}
\puthmorphism(600,0)[(\alpha \otimes \beta )\otimes \gamma `(\beta \otimes
\alpha )\otimes \gamma `\com_{\beta ,\alpha }\otimes 1_{\gamma }]{1000}{-1}a
\puthmorphism(1600,0)[\ph{(\beta \otimes \alpha )\otimes \gamma }`(\beta
\otimes \alpha )\otimes \gamma `\qss_{\beta ,\alpha }\otimes 1_{\gamma }]
{1000}{1}a
\putvmorphism(600,0)[\ph{(\alpha \otimes \beta )\otimes \gamma }`\alpha \otimes
(\beta \otimes \gamma )`\ass_{\alpha ,\beta ,\gamma }]{600}{1}l
\putvmorphism(600,-600)[\ph{\alpha \otimes (\beta \otimes \gamma )}`\alpha
\otimes (\beta \otimes \gamma )`1_{\alpha }\otimes \qss_{\beta ,\gamma }]
{600}{-1}l
\putvmorphism(600,-1200)[\ph{\alpha \otimes (\beta \otimes \gamma )}`\ph{
(\alpha \otimes \beta )\otimes \gamma }`\ass_{\alpha ,\beta ,\gamma }]
{600}{-1}l
\putvmorphism(2600,0)[\ph{(\beta \otimes \alpha )\otimes \gamma }`\beta \otimes
(\alpha \otimes \gamma )`\ass_{\beta ,\alpha ,\gamma }]{600}{1}r
\putvmorphism(2600,-600)[\ph{\beta \otimes (\alpha \otimes \gamma )}`\beta
\otimes (\alpha \otimes \gamma )`1_{\beta }\otimes \qss_{\alpha ,\gamma }]
{600}{-1}r
\putvmorphism(2600,-1200)[\ph{\beta \otimes (\alpha \otimes \gamma )}`\ph{(
\beta \otimes \alpha )\otimes \gamma }`\ass_{\beta ,\alpha ,\gamma }]{600}{-1}r
\puthmorphism(600,-1800)[(\alpha \otimes \beta )\otimes \gamma `(\alpha \otimes
\beta )\otimes \gamma `\qss_{\alpha ,\beta }\otimes 1_{\gamma }]{1000}{-1}b
\puthmorphism(1600,-1800)[\ph{(\alpha \otimes \beta )\otimes \gamma }`(\beta
\otimes \alpha )\otimes \gamma `\com_{\alpha ,\beta }\otimes 1_{\gamma }]
{1000}{1}b
\end{diagram}$
\caption{The First Decagon diagram.}
\end{figure}
\begin{figure}[h]
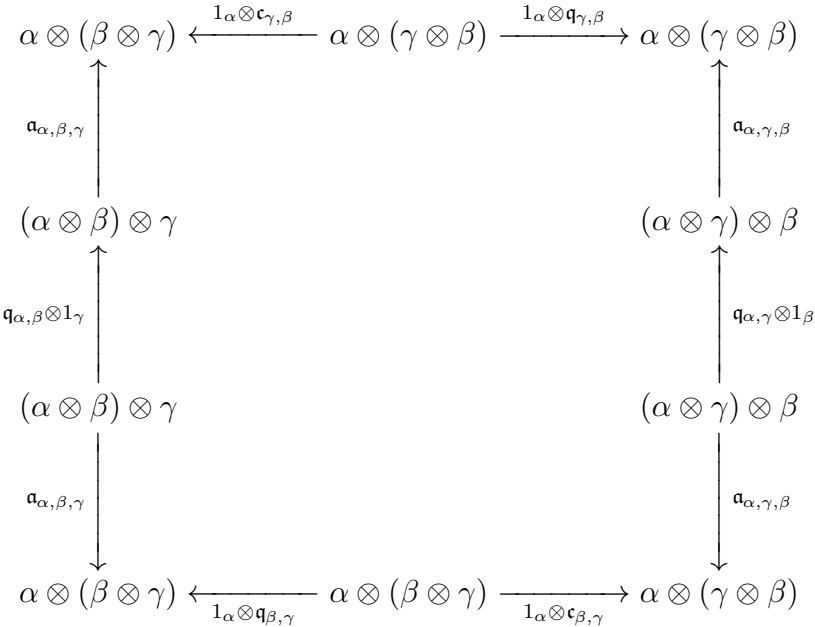

$\begin{diagram}
\puthmorphism(600,0)[\alpha \otimes (\beta \otimes \gamma )`\alpha \otimes
(\gamma \otimes \beta )`1_{\alpha }\otimes \com_{\gamma ,\beta }]{1000}{-1}a
\puthmorphism(1600,0)[\ph{\alpha \otimes (\gamma \otimes \beta )}`\alpha
\otimes (\gamma \otimes \beta )`1_{\alpha }\otimes \qss_{\gamma ,\beta}]
{1000}{1}a
\putvmorphism(600,0)[\ph{\alpha \otimes (\beta \otimes \gamma )}`(\alpha
\otimes \beta )\otimes \gamma `\ass_{\alpha ,\beta ,\gamma }]{600}{-1}l
\putvmorphism(600,-600)[\ph{(\alpha \otimes \beta )\otimes \gamma }`(\alpha
\otimes \beta )\otimes \gamma `\qss_{\alpha ,\beta }\otimes 1_{\gamma }]
{600}{-1}l
\putvmorphism(600,-1200)[\ph{(\alpha \otimes \beta )\otimes \gamma }`\ph{
\alpha \otimes (\beta \otimes \gamma )}`\ass_{\alpha ,\beta ,\gamma }]{600}{1}l
\putvmorphism(2600,0)[\ph{\alpha \otimes (\gamma \otimes \beta )}`(\alpha
\otimes \gamma )\otimes \beta `\ass_{\alpha ,\gamma ,\beta }]{600}{-1}r
\putvmorphism(2600,-600)[\ph{(\alpha \otimes \gamma )\otimes \beta }`(\alpha
\otimes \gamma )\otimes \beta `\qss_{\alpha ,\gamma }\otimes 1_{\beta }]
{600}{-1}r
\putvmorphism(2600,-1200)[\ph{(\alpha \otimes \gamma )\otimes \beta }`\ph{
\alpha \otimes (\gamma \otimes \beta )}`\ass_{\alpha ,\gamma ,\beta }]{600}{1}r
\puthmorphism(600,-1800)[\alpha \otimes (\beta \otimes \gamma )`\alpha \otimes
(\beta \otimes \gamma )`1_{\alpha }\otimes \qss_{\beta ,\gamma }]{1000}{-1}b
\puthmorphism(1600,-1800)[\ph{\alpha \otimes (\beta \otimes \gamma )}`\alpha
\otimes (\gamma \otimes \beta )`1_{\alpha }\otimes \com_{\beta ,\gamma }]
{1000}{1}b
\end{diagram}$
\caption{The Second Decagon diagram.}
\end{figure}

The underlying binary tree groupoid is the groupoid of numbered RB trees
denoted $\NRBTree $. The objects are RB trees (abbreviated NRB trees)
of any length $n$ with each leaf assigned a distinct number from $\{ 1,
...,n\} $. These numbers we write above the leaf levels in the level sequence
for the underlying RB tree. An example is given in figure 35.
\begin{figure}[h]
\hspace*{2cm}
\epsfxsize=80pt
\epsfbox{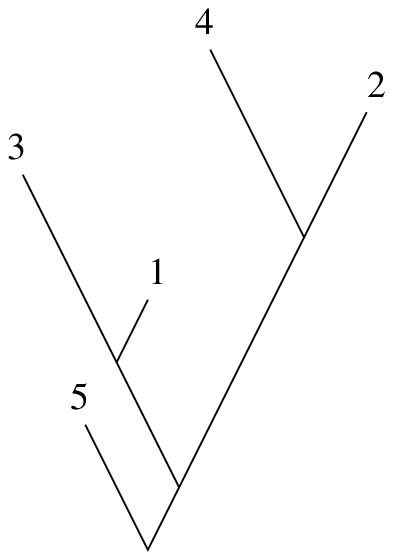}
\caption{The NRB tree given by $\COLbin{5\ph{0}3\ph{3}1\ph{1}
4\ph{5}2}{206341857}$.}
\end{figure}
The primitive arrows are inherited from $\RBTree $ togther with additional
primitive arrows corresponding to $\com $ defined as follows. Let $B$ be an
NRB tree. For any internal node $p$ let $m<p$ be the first internal node to
the left of $p$ and $n<p$ the first internal node to the right of $p$; we can
interchange the subsequences between $m$ and $p$ and $p$ and $n$. Pictorially
this corresponds to swaping the two subtrees rooted at the terminates of the
pivot $p$. Since $\com $ is symmtric the levels of the terminates is
irrelevant. The objects of $\NRBTree $ induce functors as in the $\RBTree $
situation with the difference that the objects have been permuted according to
the numbers on the leaves. If $a_{1},..,a_{n}$ is the sequence of leaf numbers
then the permutation is given by $(1,...,n)\mapsto (a_{1},...,a_{n})$.
The Square diagram and the Hexagon diagram originate from those diagrams in
$\NRBTree $ given in figures 36 and 37 respectively.
\begin{figure}[h]
\hspace*{2cm}
\epsfxsize=150pt
\epsfbox{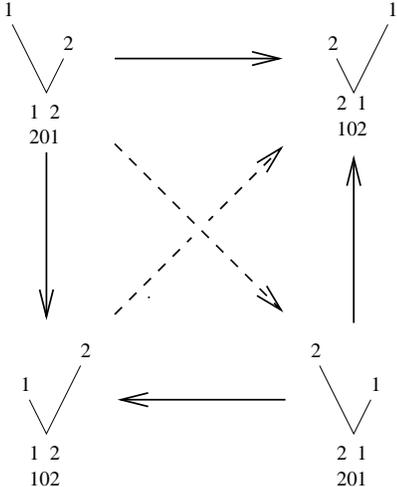}
\caption{Diagram in $\NRBTree $ underlying the Square diagram.}
\end{figure}
\begin{figure}[h]
\hspace*{2cm}
\epsfxsize=300pt
\epsfbox{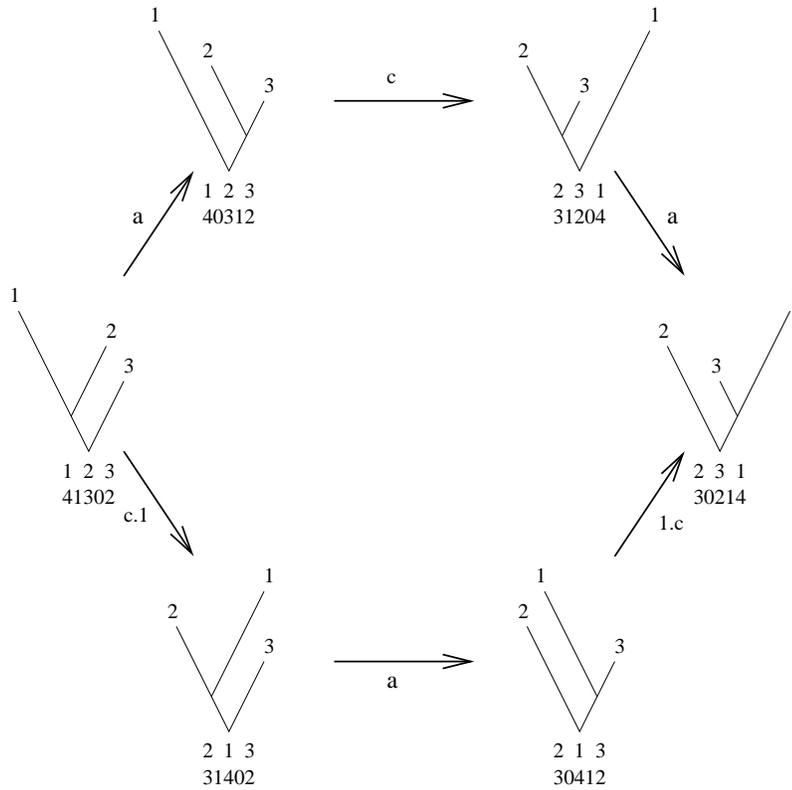}
\caption{Diagram in $\NRBTree $ underlying the Hexagon diagram.}
\end{figure}
\begin{figure}[h]
\hspace*{2cm}
\epsfxsize=250pt
\epsfbox{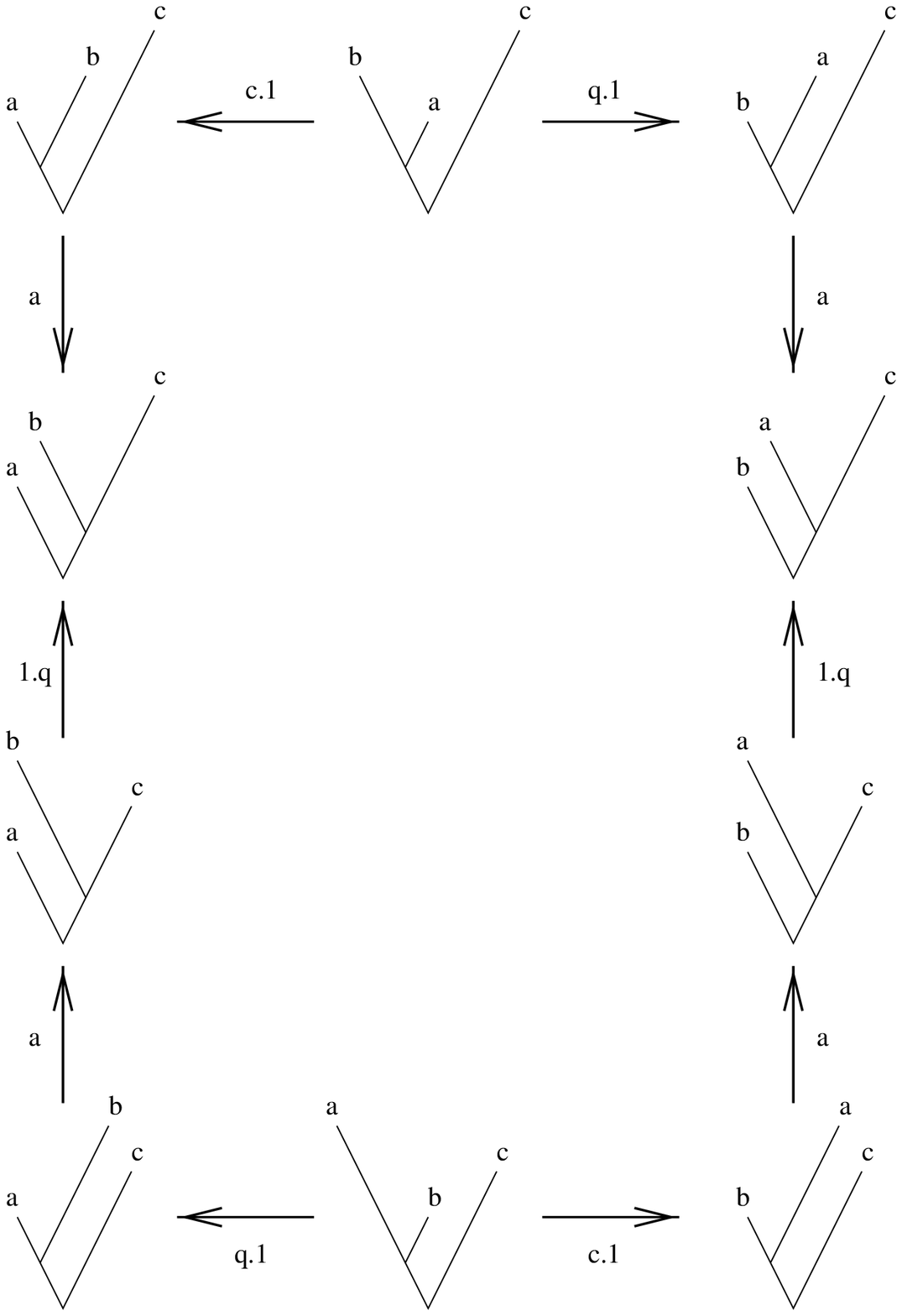}
\caption{Diagram in $\NRBTree $ underlying the First Decagon diagram.}
\end{figure}
\begin{figure}[h]
\hspace*{2cm}
\epsfxsize=250pt
\epsfbox{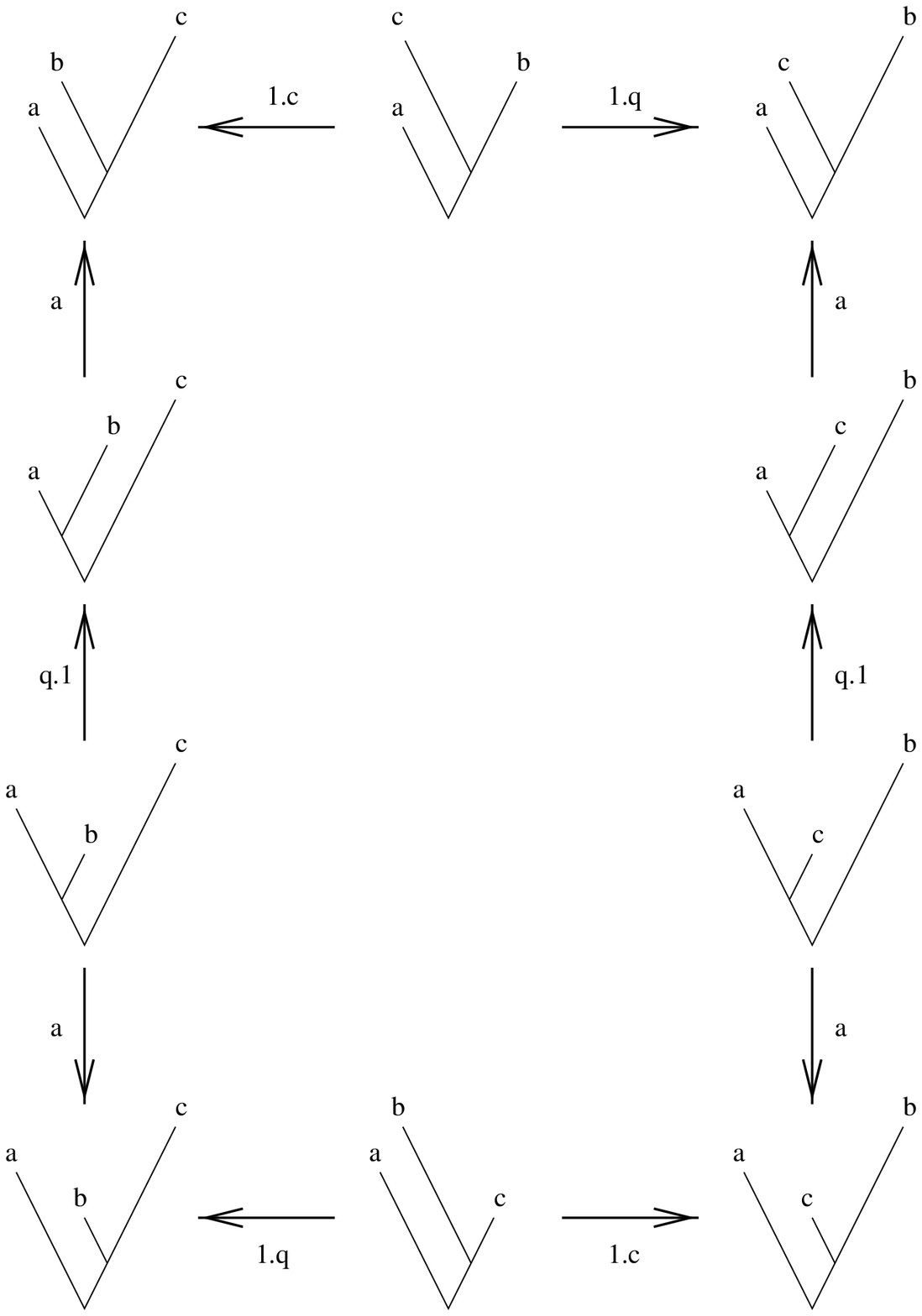}
\caption{Diagram in $\NRBTree $ underlying the Second Decagon diagram.}
\end{figure}
Note that the Hexagon diagram is independent of the terminate levels. The
diagram of figure 37 is one out of $3!$ diagrams of length three. This is why
$\com $ must be commutative or in other words indifferent to the levels
of the terminates. Furthermore, because of the Square diagram the
natural isomorphisms $\com \qss :\otimes \rightarrow \otimes \tau $ and
$\qss \com :\otimes \rightarrow \otimes \tau $ are equal, symmetric and
correspond to binary tree operations that swap terminates without
interchanging levels.

The groupoid $\NRBTree $ has the properties of the following Proposition.
\begin{prop}$\ph{nothing}$
\begin{enumerate}
\item{Given two NRB trees of length $n$ then there is a
finite sequence of primitive arrows transforming one into the other.}
\item{Every NRB tree of length $n$ is the source of at most $3(n-1)$
distinct primitive arrows.}
\item{Let $T(n)$ denote the number of RB trees of length $n$. The
number of NRB trees of length $n$ is $n!T(n)$.}
\end{enumerate}
\end{prop}
\pf\ (i) Let $B$ and $B^{\prime }$ be two NRB trees of length $n$. Every
permutation is the product of a finite number of transpositions. Hence the
result follows if every transposition of $B$ to $B^{\prime }$ can be
constructed from a sequence of primitive arrows. We construct (by Proposition
2) a sequence of internal arrows from $B$ to an NRB tree with the $i$ and
$i+1$ terminates of the branch at level $n-1$. Next we swap the terminates.
Finally we construct a sequence of internal arrows to $B^{\prime }$ (again
by Proposition 2).\\
(ii) There are a maximum of $2(n-1)$ possible distinct primitive RB arrows.
Every internal node admits a primitive arrow corresponding to a
swap. There are $n-1$ internal nodes.\\
(iii) There are $T(n)$ choices of RB trees. Each can be numbered in
$n!$ ways.

The final step is the coherence result.
\begin{thm}
  A $\qss $--pseudomonoidal category $(\cC ,\otimes ,\ass ,\qss )$ where
  $\qss $ is symmetric and $\com :\otimes \rightarrow \otimes \tau $ is
  a symmetric natural isomorphism for commutativity is symmetric $\qss
  $--pseudomonoidal coherent if and only if the Square diagram (figure
  31) and the Hexagon diagram (figure 32) and two Decagon diagrams
  (figures 33 and 34) all commute.
\end{thm}
\pf\ We first show that every iterate of $\com $ can be repalced by a sequence
of interates of $\qss $, $\ass $ and those $\com $ corresponding to
transposition of adjacent leaves. Consider an arrow $a \rightarrow b$ that is
an interate of $\com $. We define its rank $n$ to be the minimum number of
adjacency transpositions required to generate the permutation. We inductively
reduce to the desired sequence. If $n=2$ then we have our desired sequence.
Otherwise $n>2$ and we let $p$ be the pivot point in $a$ (and hence in $b$)
for the arrow $a\rightarrow b$. Let $k,l$ be the terminates of $p$. At least
one of $k,l$ is an internal node. Suppose it is $k$. If $k>l$ then we replace
the arrow using the Square diagram. Now we have $k<l$. We use the Hexagon
diagram to replace $\It_{p}(\com )$ with a sequence of five arrows where the
only iterates of $\com $ are $\It_{p}(\com .1)$ and $\It_{p}(\com )$. These
iterates have rank strictly less than $n$. Continuting this procedure
inductively we eventually replace $a\rightarrow b$ with a sequence of arrows
containing $2^{n}$ adjacent transpositions. These transpositions are the only
iterates of commutativity.

The second step of the proof is to show that two alternative sequences of
arrows, each containing exactly one adjacent transposition (iterate of $\com $)
between objects encloses a region that commutes. Let two such sequences
be $a_{1}\rightarrow \cdots \rightarrow a_{r}$ with adjacent
transposition $\It_{k}(\com ):a_{i}\rightarrow a_{i+1}$ pivoting about
$k$ with terminate levels $a,b$; and $b_{1}\rightarrow \cdots
\rightarrow b_{s}$ with adjacent transposition $\It_{l}(\com ):b_{i^{\prime }}
\rightarrow b_{i^{\prime }+1}$ pivoting about $l$ with terminate levels $c,d$.
Also $a_{1}=b_{1}$ and $a_{r}=b_{s}$. We assume that $k>l$, $c$ is to the left
of $a$ and $d$ is to the right of $b$. The other seven possibilities follow
similarly. A construction demonstrating the region enclosed commutes is given
in figure 40.
\begin{figure}[h]
\hspace*{2cm}
\epsfxsize=350pt
\epsfbox{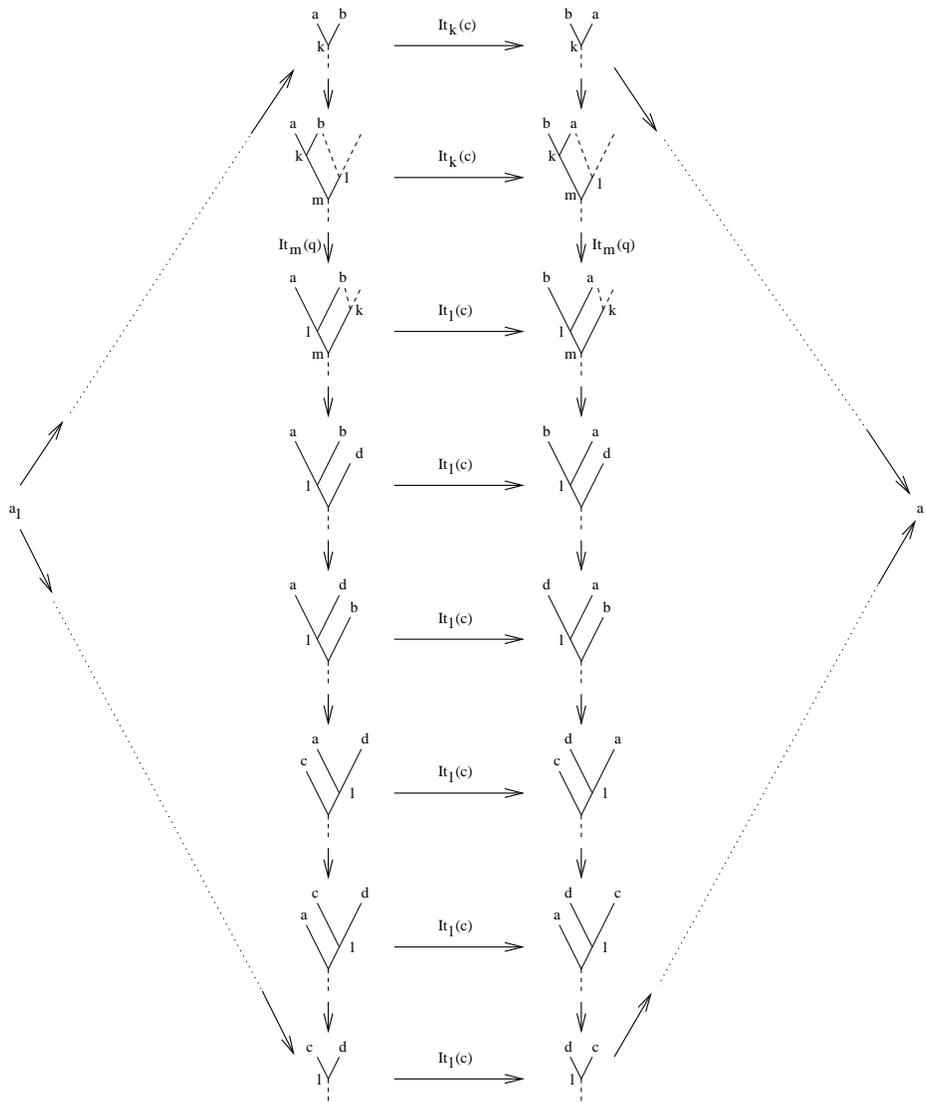}
\caption{The equivalence of two adjacent position transpositing sequences.}
\end{figure}
The sequence of arrows around the top is $a_{1}\rightarrow \cdots \rightarrow
a_{r}$. The sequence around the bottom is $b_{1}\rightarrow \cdots \rightarrow
b_{s}$. Each vertical arrow stands in for a sequence of arrows to keep the
size of the diagram acceptable. Starting from the top we proceed to describe
the ladder of regions. The sides of the top region are identical and each
arrow keeps $a$, $b$ and $k$ fixed. This region commutes by naturality and
Theorem 2. The next region down is the Square diagram. The next down has
identical sides composed of arrows keeping $a$, $b$ and $l$ fixed. This
commutes by Theorem 2 and naturality. The next region down is the first Decagon
diagram. The next region down has identical sides composed of arrows keeping
$a$, $d$ and $l$ fixed. This region commutes by Theorem 2 and naturality.
The next region down is the second Decagon diagram. The final region at the
base has identical sides composed of arrows keeping $c$, $d$ and $l$ fixed.
This region commutes by Theorem 2 and naturality. The two side
regions commute by Theorem 2.

The final step is to show that any diagram $D$ of rank $n+1$ commutes. This
proceeds in a similar way to the proof of Theorem 2 part two. We replace every
interate of $\com $ with a sequence containing only iterates of $\com $ that
are adjacent transpositions. We define the adjacent transpositions between two
objects $a$ and $b$ to be the arrows $\tau_{i}=(i\ i+1):a\rightarrow b$
swapping positions $i$ and $i+1$ given by any sequence of arrows
between the objects $a$ and $b$ with only one terminal arrow. This definition
is well--defined by the previous paragraph. Note that the adjacent
transpositions are families of arrows index by source and target. We partition
$D$ into sections where each section is a sequence of arrows containing only
one terminal arrow. Let $a_{0},...,a_{r}$ be the boundary vertices of these
sections with $a_{r+1}=a_{0}$. Each sequence $a_{k}\rightarrow \cdots
\rightarrow a_{k+1}$ corresponds to an adjacent transposition $\tau_{i}:a_{k}
\rightarrow a_{k+1}$. We now show that the adjacent transpositions satisfy the
generating relations of the permutation group $S_{n+1}$.
\begin{eqnarray}
\tau_{i}^{2}=1~, & & i=1,...,n \\
\tau_{i}\tau_{j}=\tau_{j}\tau_{i}~, & & 1\leq i<j\leq n\\
\tau_{i}\tau_{i+1}\tau_{i}=\tau_{i+1}\tau_{i}\tau_{i+1}~, & \ph{space} &
i=1,...,n-1
\end{eqnarray}
Properties (i), (ii) and (iii) are proved by an obvious modification to the
verification in the proof of Theorem 2 part two. For property (iii) we note
that the Dodecagon diagram  embedded in figure 21 is replace by the diagram of
figure 41.
\begin{figure}[h]
\epsfxsize=450pt
\epsfbox{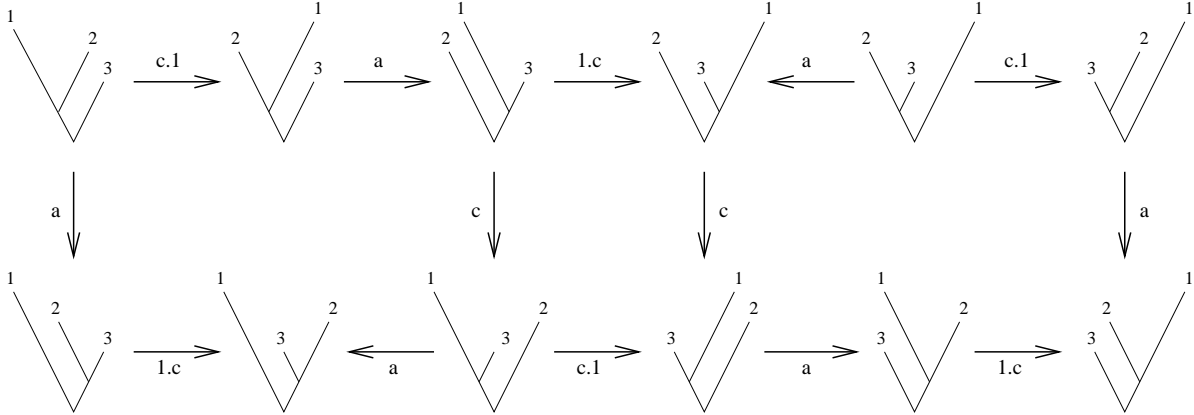}
\caption{Verification of the property $\tau_{i}\tau_{i+1}\tau_{i}=\tau_{i+1}
\tau_{i}\tau_{i+1}$.}
\end{figure}
The middle region is a naturality square and the two side regions are Hexagon
diagrams. This completes the proof.

Symmetric Pseudomonoidal coherence follows as a corollary to Theorem 5. We give
the definition leaving the details to the reader.
\begin{defn}
  A symmetric pseudomonoidal category is a quadruple $(\cC ,\otimes
  ,\ass ,\com )$ where $(\cC ,\otimes ,\ass )$ is a pseudomonoidal
  category, $\qss $ and $\com $ are symmetric, the Square diagram,
  Hexagon diagram and Decagon diagrams (extended to branch points) hold
  and the diagrams of figure 42 hold.
\end{defn}
\begin{figure}[h]
\hspace*{2cm}
$\begin{diagram}
\putsquare<1`1`1`1;1200`600>(600,0)[(\alpha \otimes \beta )\otimes (\gamma
\otimes \delta )`(\beta \otimes \alpha )\otimes (\gamma \otimes \delta )`
(\alpha \otimes \beta )\otimes (\gamma \otimes \delta )`(\beta \otimes \alpha )
\otimes (\gamma \otimes \delta );\com_{\alpha ,\beta }\otimes 1_{\gamma \otimes
\delta }`\qss_{\alpha ,\beta ,\gamma ,\delta }`\qss_{\beta ,\alpha ,\gamma
,\delta }`\com_{\alpha ,\beta }\otimes 1_{\gamma \otimes \delta }]
\putsquare<1`1`1`1;1200`600>(600,-1000)[(\alpha \otimes \beta )\otimes (\gamma
\otimes \delta )`(\alpha \otimes \beta )\otimes (\delta \otimes \gamma )`
(\alpha \otimes \beta )\otimes (\gamma \otimes \delta )`(\alpha \otimes \beta )
\otimes (\delta \otimes \gamma );1_{\alpha \otimes \beta }\otimes \com_{\gamma
,\delta }`\qss_{\alpha ,\beta ,\gamma ,\delta }`\qss_{\alpha ,\beta ,\delta
,\gamma }`1_{\alpha \otimes \beta }\otimes \com_{\gamma ,\delta }]
\end{diagram}$
\caption{Additional square diagrams required for symmetric pseudomonoidal
coherence.}
\end{figure}

\section{Symmetric $\qss $--Monoidal Categories}                 %

The difficulty with unital pseudomonoidal structures is that the notion
of identity is handled monoidally. This rather defeats the underlying
motivation for considering pseudomonoidal categories in the first place.
Including the notion of identity has proven to be a delicate balance. A
pseudomonoidal category is too general to incorporate this notion. On the
other hand the $\qss $--pseudomonoidal category has too many conditions
resulting in the severe conditions that $\qss $ is symmetry,
$\qss_{\alpha ,\beta }\otimes 1_{\gamma } =\qss_{\alpha \otimes \beta
  ,\gamma }$ and $1_{\alpha }\otimes \qss_{\beta ,\gamma }=\qss_{\alpha
  ,\beta \otimes \gamma }$ for all objects $\alpha ,\beta ,\gamma $.
Moreover, these conditions are not represented by any binary tree
diagram in the coherence groupoid. The right level of structure is that
of a $\qss $--braided premonoidal category. Well not quite. The
existence of a commutativity natural isomorphism is also required. In
the absence of commutativity one is faced with requiring an infinite
number of diagrams to hold in order to guarantee coherence.
\begin{defn}
A symmetric $\qss $--monoidal category is an octuple $(\cC ,\otimes ,\ass ,
\qss ,\com ,\lid ,\rid ,e)$ where $(\cC ,\otimes ,\ass  ,\qss )$ is a
$\qss $--braided pseudomonoidal category, $(\cC ,\otimes ,\ass ,\com )$ is a
symmetric pseudomonoidal category, $e$ is an object of $\cC $
called the identity and $\lid :e\otimes \ubar \rightarrow 1$ and $\rid :\ubar
\otimes \rightarrow 1$ are natural isomorphisms satisfying the small and large
$\qss $--Triangle diagrams (figures 43 and 44)
\end{defn}
\begin{figure}[h]
$\begin{diagram}
\puthmorphism(0,0)[(\alpha .e).\beta `\alpha .(e.\beta )`\ass_{\alpha ,e,
\beta }]{700}{1}a
\puthmorphism(700,0)[\ph{\alpha .(e.\beta )}`\alpha .(e.\beta )`1_{\alpha }.
\qss_{e,\beta }]{700}{-1}a
\puthmorphism(1400,0)[\ph{\alpha .(e.\beta )}`(\alpha .e).\beta `\ass_{\alpha
,e,\beta }]{700}{-1}a
\puthmorphism(2100,0)[\ph{(\alpha .e).\beta }`(\alpha .e).\beta `\qss_{\alpha
,e}.1_{\beta }]{700}{-1}a
\puthmorphism(2800,0)[\ph{(\alpha .e).\beta }`\alpha .(e.\beta )`\ass_{\alpha
,e,\beta }]{700}{1}a
\putvmorphism(0,0)[\ph{(\alpha .e).\beta }`\ph{\alpha .\beta }`\rid_{\alpha }
.1_{\beta }]{700}{1}l
\putvmorphism(3500,0)[\ph{\alpha .(e.\beta )}`\ph{\alpha .\beta }`1_{\alpha }.
\lid_{\beta }]{700}{1}r
\puthmorphism(0,-700)[\alpha .\beta `\alpha .\beta `\qss_{\alpha ,\beta }]
{3500}{-1}b
\end{diagram}$
\caption{The Large $\qss $--Triangle diagram.}
\end{figure}
\begin{figure}[h]
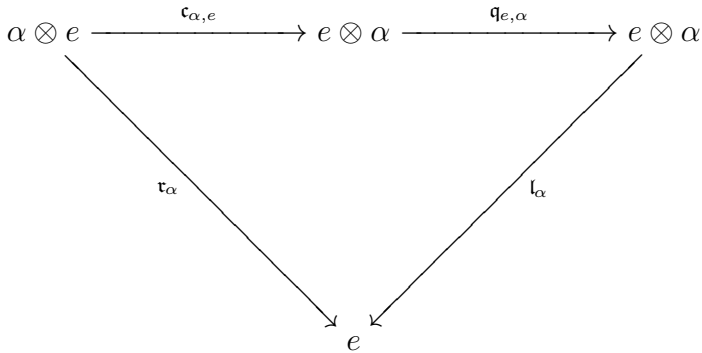

$\begin{diagram}
\puthmorphism(500,0)[\ph{\alpha \otimes e}`\ph{e\otimes \alpha }`
\com_{\alpha ,e}]{1000}{1}a
\puthmorphism(1500,0)[e\otimes \alpha `e\otimes \alpha
`\qss_{e,\alpha }]{1000}{1}a
\putVtriangle<0`1`1;1000>(500,-1000)[\alpha \otimes e`\ph{\alpha \otimes e}`e;
`\rid_{\alpha }`\lid_{\alpha }]\end{diagram}$
\caption{The Small $\qss $--Triangle diagram.}
\end{figure}
The $\qss $--Triangle diagrams collapse to the Triangle diagrams (of monoidal
categories) when $\qss =1$. There are no $\qss $--Triangle diagrams
correspondong to the redundant triangle diagrams that were originally in the
definition of a monoidal category. The redundancy amongst the original Triangle
diagrams was pointed out by Kelly \cite{gk}.

The underlying binary tree category is the groupoid of numbered RRB trees with
nodules, denoted $\NRRNBTree $. As we did earlier we represent a nodule by
attaching a small circle to the leaf. An RRB tree with nodules is represnted by
an ordered pair of linear orderings. The first entry gives the branch
structure. The second entry the leaf structure. Leafs with nodules are
represented by placing a line under its level. The length of an RRB tree with
nodules is the number of leaves less the number of nodules. These trees under
the $\can $ functor give functors where the nodules stand in for the identity.
A numbered RRNB tree of length $n$ attaches a number $1,...,n$ uniquely to
each of the nodule free leaves.

The arrows of $\NRRNBTree $ are generated by the primitive arrows inherited
from $\NRRBTree $ together with primitive arrows corresponding to left and
right identity that we now define. Given an RRB tree $B$ with nodules having
$n$ leaves; we can prune a nodule at level $n+1$ that is the left terminate of
a branch at level $n-1$ whose right terminate leaf is at level $n$. This
leaves the leaf at level $n$. All other leaf levels are lowered by two. We can
prune a nodule at level $n$ that is the right terminate of a branch at level
$n-1$ whose left terminate leaf is at level $n+1$. This leaves the leaf at
level $n$. All other leaf levels are lowered by two. The grafting arrows
are given by the converse of the prunning operation described. These primitive
arrows are given in figure 45.
\begin{figure}[h]
\hspace*{2cm}
\epsfxsize=200pt
\epsfbox{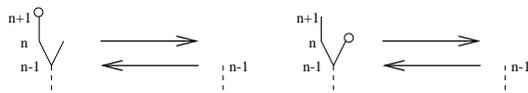}
\caption{The prunning and grafting primitive arrows of $\RRNBTree $.}
\end{figure}
The diagrams in $\RRNBTree $ underlying the $\qss $--Triangle diagrams are
given by figures 46 and 47.
\begin{figure}[h]
\hspace*{2cm}
\epsfxsize=300pt
\epsfbox{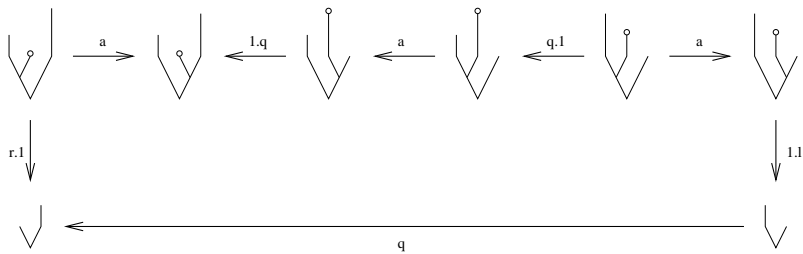}
\caption{Diagram in $\NRRNBTree $ underlying the Large $\qss $--Triangle
diagram.}
\end{figure}
\begin{figure}[h]
\hspace*{2cm}
\epsfxsize=100pt
\epsfbox{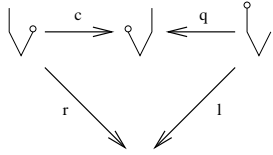}
\caption{Diagram in $\NRRNBTree $ underlying the Small $\qss $--Triangle
diagram.}
\end{figure}

\begin{prop}$\ph{nothing}$
\begin{enumerate}
\item{Given two numbered RRB trees with nodules of length $n$ then there is a
finite sequence of primitive arrows transforming one into the other.}
\item{Every numbered RRB tree with $m$ nodules of length $n$ is the source of
at most $2(n+m)$ distinct primitive arrows.}
\item{The number of RRB trees with $m$ nodules of length $n$ is
$2(n+m)!(n+m-1)!n!$.}
\end{enumerate}
\end{prop}
\pf\ (i) Let $B$ and $B^{\prime }$ be any two numbered RRB trees with nodules
of length $n$. Choosing each nodule of $B$ in turn we can rearrange using a
sequence of $\NRRBTree$ primitive arrows so that it is a terminate of the
highest branch and this branch's terminates are on the next two levels up.
Clearly we can prune the nodule. Applying the same procedure to $B^{\prime }$
we can construct a similar sequence of primitive arrows. Clearly $B$ and
$B^{\prime }$ are connected by a sequence of primitive arrows (Proposition 3).
\\
(ii) There are at most $n+m-1$ $\RRBTree $ primitive arrows. The only
possiblities for prunning or grafting is when the highest branch has a
terminate at the lowest level. There are at most two possibilities in this
situtation. Finally every branch, of which there are $n+m-1$, admits a
reflection.\\
(iii) There are $(n+m)1(n+m-1)!$ RRB trees of length $n+m$. There are $n!$
ways of numbering the nodule free leafs. Nodules may only be grafted to the
leaf at level $n+m$. There are precise two ways of doing this.

Now that we have a groupoid of binary trees describing the diagrams of a
symmetric $\qss $--monoidal category we formulate the notion of coherence in
the by now standard way.
\begin{thm}
A symmetric $\qss $--braided pseudomonoidal category $(\cC ,\otimes ,\ass ,\qss
,\com )$ with an object $e$ and natural isomorphisms $\lid :e\otimes \ubar
\rightarrow 1$ and $\rid :\ubar \otimes e\rightarrow 1$ is symmetric
$\qss $--monoidal coherent if and only if the $\qss $--Triangle diagrams
(figures 43 and 44) both
commute.
\end{thm}
\pf\ The proof is by induction. Define the rank of a vertex to be the number of
leafs and the rank of a diagram to be the maximum of its vertex ranks.
The theorem holds for rank $3$ diagrams. Suppose all diagrams with rank $n$
are coherent. Let $D$ be a diagram of rank $n+1$ with vertices
$a_{0},a_{1},...,a_{r}$ reading around the outside and put $a_{r+1}=a_{0}$. If
every vertex of $D$ has rank $n+1$ then there are no primitive arrows
prunning/grafting nodules. Hence $D$ commutes by Theorem 5. Otherwise we
divide $D$ into maximal sequences where all the vertices have rank $n+1$
alternating with none of the vertices have rank $n+1$. Let one such sequence
be where all the vertices have rank $n+1$ be $a_{k}\rightarrow \cdots
\rightarrow a_{l}$. The arrows $a_{k-1}\rightarrow a_{k}$, $a_{l}\rightarrow
a_{l+1}$ fall into three cases.\\
\textbf{Case (a):} The arrows do not prune/graft a nodule into the first or
second position. If both arrows are of this type then we substitute the
sequence $a_{k-1}\rightarrow \cdots \rightarrow a_{l+1}$ according to figure
48.
\begin{figure}[h]
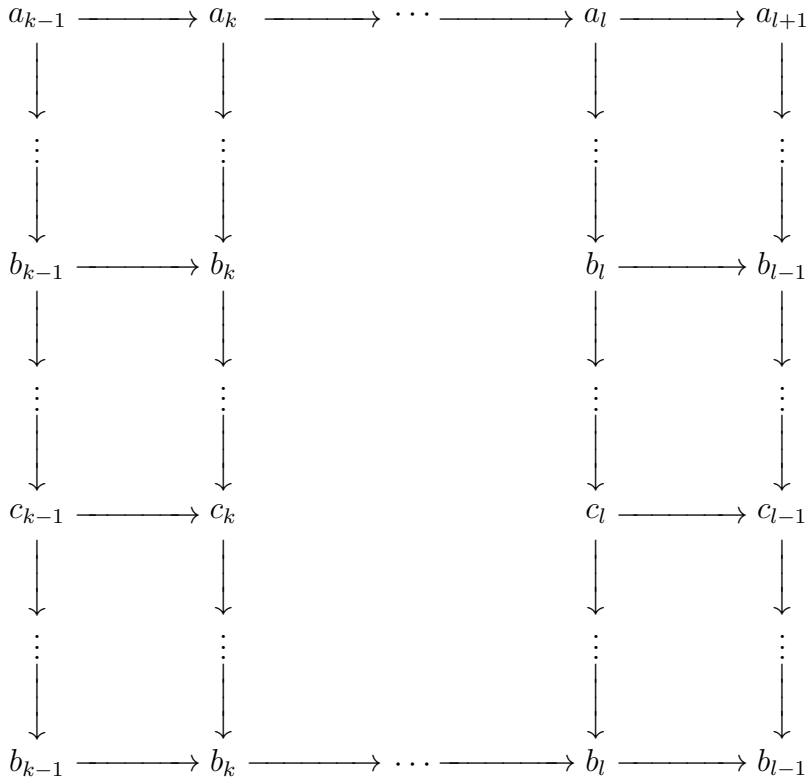

$\begin{diagram}
\putsquare<1`1`1`0;600`400>(0,0)[a_{k-1}`a_{k}`\vdots `\vdots ;```]
\puthmorphism(600,400)[\ph{a_{k-1}}`\cdots `]{600}{1}a
\puthmorphism(1200,400)[\ph{\cdots }`\ph{a_{l}}`]{600}{1}a
\putsquare<1`1`1`0;600`400>(1800,0)[a_{l}`a_{l+1}`\vdots `\vdots ;```]
\putsquare<0`1`1`1;600`400>(0,-400)[\ph{\vdots }`\ph{\vdots }`b_{k-1}`b_{k};
```]
\putsquare<0`1`1`1;600`400>(1800,-400)[\ph{\vdots }`\ph{\vdots }`b_{l}`b_{l-1};
```]
\putsquare<0`1`1`0;600`400>(0,-800)[\ph{b_{k-1}}`\ph{b_{k}}`\vdots `\vdots ;
```]
\putsquare<0`1`1`0;600`400>(1800,-800)[\ph{b_{l}}`\ph{b_{l+1}}`\vdots `\vdots ;
```]
\putsquare<0`1`1`1;600`400>(0,-1200)[\ph{\vdots }`\ph{\vdots }`c_{k-1}`c_{k};
```]
\putsquare<0`1`1`1;600`400>(1800,-1200)[\ph{\vdots }`\ph{\vdots }`c_{l}`
c_{l-1};```]
\putsquare<0`1`1`0;600`400>(0,-1600)[\ph{c_{k-1}}`\ph{c_{k}}`\vdots `\vdots ;
```]
\putsquare<0`1`1`0;600`400>(1800,-1600)[\ph{c_{l}}`\ph{c_{l+1}}`\vdots `\vdots
;```]
\putsquare<0`1`1`1;600`400>(0,-2000)[\ph{\vdots }`\ph{\vdots }`b_{k-1}`b_{k};
```]
\puthmorphism(600,-2000)[\ph{b_{k}}`\cdots `]{600}{1}b
\puthmorphism(1200,-2000)[\ph{\cdots }`\ph{b_{l}}`]{600}{1}b
\putsquare<0`1`1`1;600`400>(1800,-2000)[\ph{\vdots }`\ph{\vdots }`b_{l}`
b_{l-1};```]
\end{diagram}$
\caption{Removal of vertices with rank $n+1$ from $D$.}
\end{figure}
The top horizontal line is the sequence $a_{k-1}\rightarrow \cdots \rightarrow
a_{l+1}$. Suppose that the highest leaf level is to the left of the nodule
being grafted by the arrow $a_{k-1}\rightarrow a_{k}$ and that this arrow
corresponds to an interate of the right identity. We construct identical
sequences of $\RRBTree $ arrows from $a_{k-1}$ to $b_{k-1}$ and from $a_{k}$
to $b_{k}$. The arrow $b_{k-1}\rightarrow b_{k}$ is an iterate of $\rid $ and
the regin enclose commutes by naturality. The vertices $b_{k-1}$ and $b_{k}$
are of the form given in figure 49.
\begin{figure}[h]
\hspace*{2cm}
\epsfxsize=200pt
\epsfbox{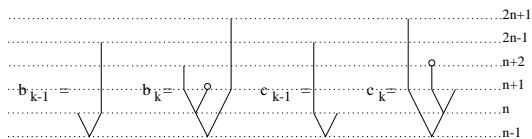}
\caption{The vertices $b_{k-1}$, $b_{k}$, $c_{k-1}$ and $c_{k}$.}
\end{figure}
The top left hand region in figure 48 commutes by naturality and Theorem 5.
The next region down is taken to be the large $\qss $--Triangle diagram where
the vertices $c_{k-1}$ and $c_{k}$ are as given in figure 49. We construct
identical sequences of $\RRBTree $ arrows from $c_{k-1}$ to $d_{k-1}$ and
from $c_{k}$ to $d_{k}$. The arrow $d_{k-1}\rightarrow d_{k}$ is an iterate of
$\rid $ where the vertices $d_{k-1}$ and $d_{k}$ are of the form given in
figure 50.
\begin{figure}[h]
\hspace*{2cm}
\epsfxsize=200pt
\epsfbox{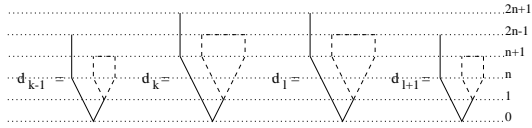}
\caption{The vertices $d_{k-1}$, $d_{k}$, $d_{l}$ and $d_{l+1}$.}
\end{figure}
The region these sequences (in figure 48) enclose commutes. If the arrow
$a_{k-1}\rightarrow a_{k}$ was an iterate of the left identity then the
above construction follows similarly with the middle left hand region of
figure 48 being a combination of the small and large $\qss $--Triangle
diagrams. The only other possiblilty for $a_{k-1}\rightarrow a_{k}$ is if the
highest leaf level is to the right of the nodule being grafted. In this case
the construction proceeds directly from $a_{k-1}$ to $d_{k-1}$ and from $a_{k}$
to $d_{k}$.

The left hand side vertical sequences of arrows in figure 46 is constructed
as for the right hand side with all the enclosed regions commuting. Finally
we connect $d_{k}$ to $d_{l}$ using a sequence of $\NRRBTree $ arrows that
do not pivot about the root. The region enclosed contains no arrows for
grafting nor prunning nodules and by Theorem 5 commutes. The sequence
running arround the bottom of the diagram from $a_{k-1}$ to $d_{k-1}$ to
$d_{l+1}$ to $a_{l+1}$ is substituted for the maximal sequence.\\
\textbf{Case (b):} The arrows do not prune/graft a nodule into the first
position but at least one of these does into the second position. Suppose
$a_{k-1}\rightarrow a_{k}$ grafts a nodule into the second position. This
arrow is one of the top horizontal arrows given in the diagrams of figure 51.
\begin{figure}[h]
\hspace*{2cm}
\epsfxsize=250pt
\epsfbox{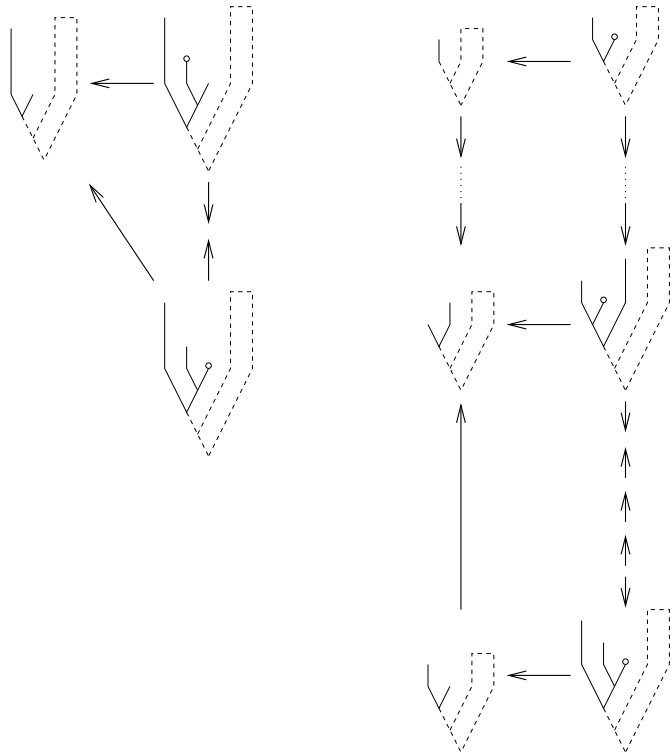}
\caption{Diagrams used to reduce from case (b) to case (a).}
\end{figure}
If it is the first diagram then we see that the small $\qss $--Triangle
diagram allows us to substitute it for a sequence grafting a nodule into the
third position. If it is the second diagram then we construct the two sides of
the top region using two identical sequences of $\RRBTree $ arrows that
do not pivot about the highest branch and maintain the first and second
postions as terminates. As before the region enclosed commutes. Finally the
bottom region is the large $\qss $--Triangle diagram. We use the sequence
running around the bottom of the diagram (of figure 51) to substitute for
$a_{k-1}\rightarrow a_{k}$ for which the primitive arrow for grafting
a nodule into the third position. We apply a similar substitution to
$a_{l}\rightarrow a_{l+1}$ if it prunes a nodule from the second position. The
new maximal sequence falls into case (b).\\
\textbf{Case (c):} At least one arrow prunes/grafts a nodule into the first
position. We replace all arrows grafting/prunning a nodule into the first with
one grafting/prunning into the second position using the small
$\qss $--Triangle diagram. The new maximal sequence falls into case (a) or (b).

We divide $D$ up into maximal sequences as before. Let $a_{k}\rightarrow \cdots
\rightarrow a_{l}$ and $a_{k^{\prime }}\rightarrow \cdots \rightarrow
a_{l^{\prime }}$ be consecutive maximal sequences. The vertices
$a_{k-1},...,a_{l+1},a_{k^{\prime }-1},...,a_{l^{\prime }+1}$ all have the
highest leaf as the left terminate of the root. We replace
the joining sequence $a_{l+1}\rightarrow \cdots \rightarrow a_{k^{\prime }-1}$
by an alternative sequence of $\NRRBTree $ primitive arrows with the
highest leaf as the left terminate of the root. The region enclosed commutes
by Theorem 5. Conituning this substitution inductively we construct a
diagram where every vertex has the left terminate of the root as the highest
leaf. Moreover, this diagram commutes if and only if $D$ commutes. Hence by the
induction hypothesis $D$ commutes. This completes the proof.

\section{Summary}                                                %

The prototype of all structures considered in this paper is that of a
premonoidal/pseudomonoidal category. This gave rise to the notion of a
$\qss $ natural automorphism which in the coherence groupoid corresponds
to interchanging the level of two branch points. Extending this to
include leaves leads to the notion of $\qss $--pseudomonoidal category.
Prohibiting the interchange of a leaf level with a branch level leads to
the notion of a $\qss $--braided pseudomonoidal category. The hierarchy of
these structures is given in table 3 together with their coherence
groupoid.
\begin{table}[h]
\hspace*{4cm}
\begin{tabular}{|c|c|}
\hline
structure & coherence groupoid \\
\hline
pseudomonoidal & $\IRBTree $\\
$\qss $--braided pseudomonoidal & $\RRBTree $\\
$\qss $--pseudomonoidal & $\RBTree $\\
monoidal & $\BTree $\\
\hline
\end{tabular}
\caption{Coherence groupoids underlying natural associativity structures.}
\end{table}
Each structure admits restricted and symmetric versions. Only the symmetric
$\qss $--symmetric pseudomonoidal category admits a $\qss $--monoidal structure
which we call symmetric $\qss $--monoidal. The relationship between all the
structures studied in this paper is depicted in figure 52.
\begin{figure}[h]
\hspace*{2cm}
\epsfxsize=350pt
\epsfbox{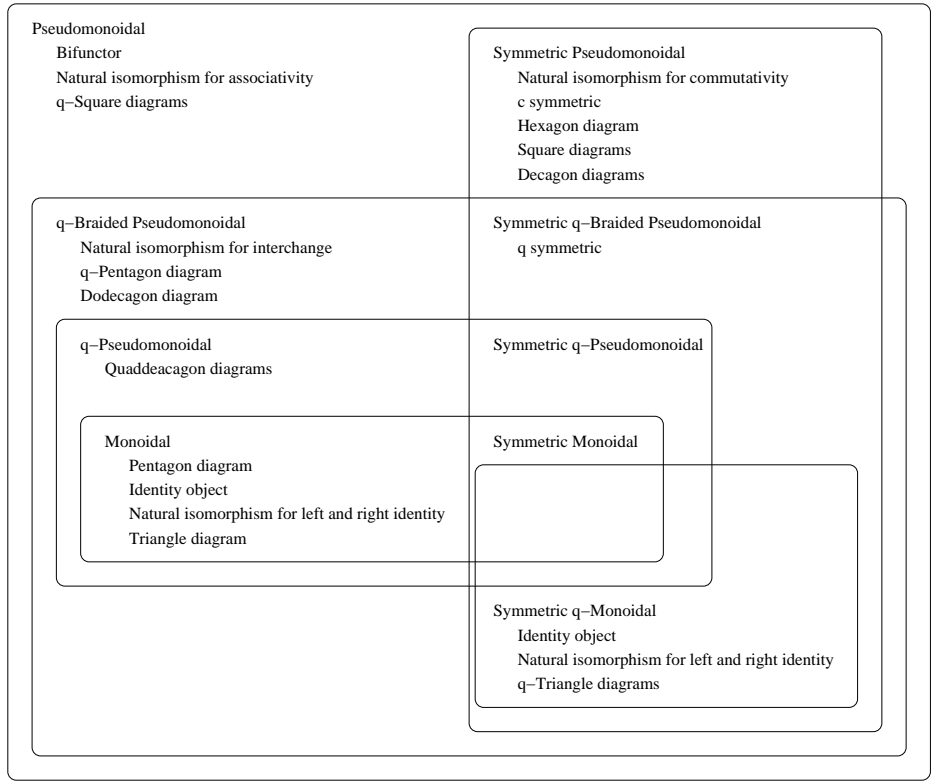}
\caption{Overview of pseudomonoidal structures.}
\end{figure}
Finally further work on braided premonoidal coherence and the breaking
of unital structure is soon to appear \cite{wj3}.

\section*{Appendix A}                                            %

A $\qss $--pseudomonoidal category allows one to ponder the existence of a
$\qss $--associahedra as a deformation of the associahedra for monoidal
categories. The $\qss $--associahedron for words of length four is the
$\qss $--Pentagon diagram. For words of length five it is the planar
diagram given in figure 53.
\begin{figure}[h]
\hspace*{2cm}
\epsfxsize=350pt
\epsfbox{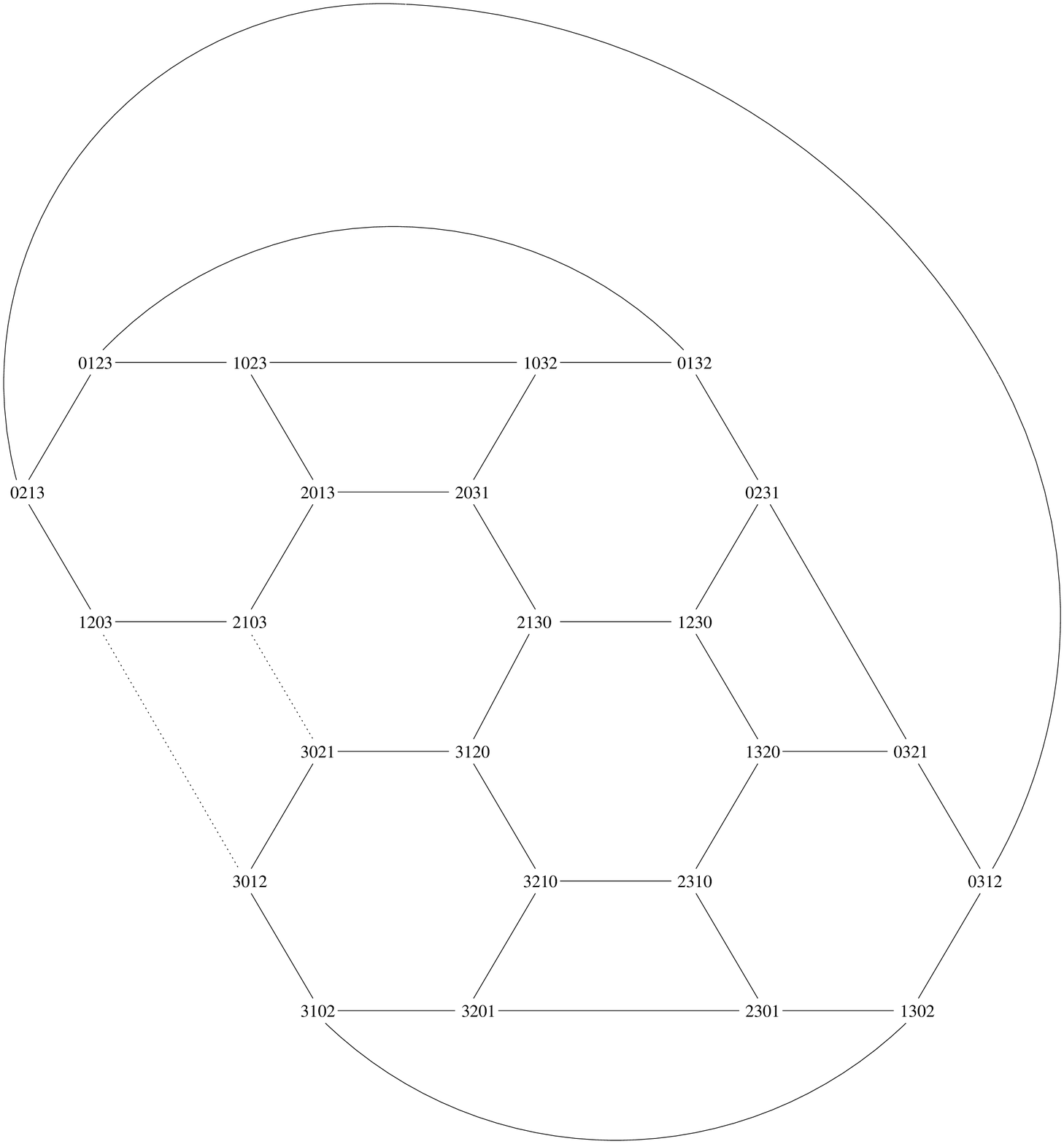}
\caption{The primitive arrows for words of length five.}
\end{figure}
Note that the dotted lines represent disallowed primitive operations
for level interchange. The hexagons correspond to the $\qss $--diagram and the
quadralaterals are natural squares. Hence the twelve vertex diagram around the
perimeter commutes. However, the diagram does not fold into a
polyhedron. The diagram wraps around a truncated polyhedron as given by
figure 54.
\begin{figure}[h]
\hspace*{2cm}
\epsfxsize=150pt
\epsfbox{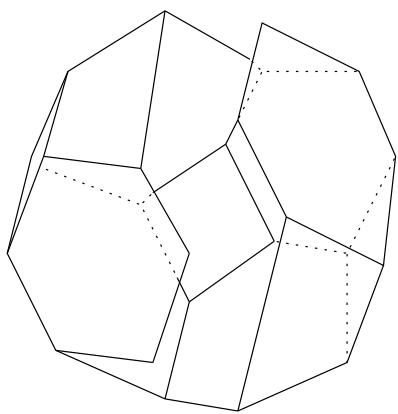}
\caption{Folding of the planar diagram of figure 51.}
\end{figure}
The missing strip corresponds to the twelve vertex diagram encircling
figure 53. This partially formed polyhedron is not what one would like
to take as the $\qss $--associahedron. Instead we utilise the
connection with permutations described in section 4. The $\qss
$--associahedron is replaced by a permutahedron. The edges to not
correspond to primitive arrows. They are, however, constructed from
primitive arrows and correspond to adjacent transpositions.

\section*{Acknowledgements}                                      %

The author wishes to thank Prof. Ross Street for supporting this
research at Macquarie University, Division of Communication and
Information Sciences, N.S.W. 2109, Australia. I would like to
acknowledge Daniel Steffen and Dr. Mark Weber for pointing out the links
between $T(n)$, tangent numbers and up/down permutations. Also I thank
Dr. Michael Batanin for some useful discussions, and Prof. Jim Stasheff
for useful suggestions on an early draft.

\section*{References}                                            %

\end{document}